\documentclass[3p,a4paper]{elsarticle}
\usepackage[utf8]{inputenc}
\usepackage{subfigure}

\title{Scaled-cPIKANs: Domain Scaling in Chebyshev-based Physics-informed Kolmogorov-Arnold Networks}

\author[myUaddress]{Farinaz Mostajeran}
\author[myUaddress]{Salah A Faroughi \corref{mycorrespondingauthor}}

\cortext[mycorrespondingauthor]{Corresponding author: salah.faroughi@utah.edu}

\address[myUaddress]{Department of Chemical Engineering, University of Utah, Salt Lake City, Utah  84112, USA
}

\date{\today}

\usepackage{mwe}

\usepackage{setspace} 
\usepackage{tabularx}
\usepackage{lineno,hyperref}
\usepackage{amsmath}
\usepackage{bm}
\usepackage{amssymb}
\usepackage[]{amsmath}
\usepackage{upgreek}
\usepackage{float}

\let\today\relax
\makeatletter
\def\ps@pprintTitle{%
    \let\@oddhead\@empty
    \let\@evenhead\@empty
    \def\@oddfoot{\footnotesize\itshape
         {Submitted preprint — January 2025} \hfill\today}%
    \let\@evenfoot\@oddfoot
    }
\makeatother

\usepackage{pgfplots}
\pgfplotsset{compat=1.5}
\usepgfplotslibrary{units}
\usetikzlibrary{spy}
\usetikzlibrary{pgfplots.groupplots}

\usepackage{epstopdf}
\usepackage{units}
\usepackage{lscape}
\usepackage{booktabs}
\usepackage{gensymb}
\usepackage{multirow}

\usepackage{graphicx}
\usepackage{caption}
\usepackage{subcaption}

\usepackage{booktabs}
\usepackage{multirow}
\usepackage{caption}
\usepackage{subcaption}

\usepackage{soul,color}
\soulregister\cite7
\soulregister\ref7
\soulregister\pageref7

\usepackage[flushleft]{threeparttable}
\usepackage{tikz}
\usepackage{placeins}
\usepackage{color}
\usepackage{hyperref}

\usepackage{stackengine}

\usepackage{array}
\usepackage{tabularx}
\usepackage{pdfpages}
\usepackage{natbib}
\usepackage{textcomp, gensymb}
\usepackage{booktabs} % Ensure this package is included in your preamble
\begin{document}

\begin{abstract}
Partial Differential Equations (PDEs) are integral to modeling many scientific and engineering problems. Physics-informed Neural Networks (PINNs) have emerged as promising tools for solving PDEs by embedding governing equations into the neural network loss function. However, when dealing with PDEs characterized by strong oscillatory dynamics over large computational domains, PINNs based on Multilayer Perceptrons (MLPs) often exhibit poor convergence and reduced accuracy.
To address these challenges, this paper introduces Scaled-cPIKAN, a physics-informed architecture rooted in Kolmogorov-Arnold Networks (KANs). Scaled-cPIKAN integrates Chebyshev polynomial representations with a domain scaling approach that transforms spatial variables in PDEs into the standardized domain \([-1,1]^d\), as intrinsically required by Chebyshev polynomials. By combining the flexibility of Chebyshev-based KANs (cKANs) with the physics-driven principles of PINNs, and the spatial domain transformation, Scaled-cPIKAN enables efficient representation of oscillatory dynamics across extended spatial domains while improving computational performance. We demonstrate Scaled-cPIKAN efficacy using four benchmark problems: the diffusion equation, the Helmholtz equation, the Allen-Cahn equation, as well as both forward and inverse formulations of the reaction-diffusion equation (with and without noisy data).
Our results show that Scaled-cPIKAN significantly outperforms existing methods in all test cases. In particular, it achieves several orders of magnitude higher accuracy and faster convergence rate, making it a highly efficient tool for approximating  PDE solutions that feature oscillatory behavior over large spatial domains. 
\end{abstract}

\begin{keyword}
    Physics-informed Neural Networks\sep%
    Kolmogorov-Arnold Network \sep% 
    Chebyshev Polynomials \sep% 
    Domain-scaling \sep%
     Variable-scaling \sep%
    Scientific Machine Learning  
\end{keyword}

\maketitle

\section{Introduction}\label{sec:Intro}

Physics-informed neural networks (PINNs) are a class of deep learning frameworks that directly integrate physical laws, often expressed by partial differential equations (PDEs), directly into the neural network training process \cite{raissi2019physics}. 
PINNs use the universal approximation capabilities of Multilayer Perceptrons (MLPs) to resolve forward and inverse problems across a wide range of domains, such as fluid mechanics \cite{cai2021physics, faroughi2023physics}, solid mechanics \cite{hu2024physics, rezaei2024finite},  biomedical engineering  \cite{bhargava2024enhancing, caforio2024physics},  materials sciences and engineering \cite{caforio2024physics,datta2022physics}, and many more \cite{toscano2024pinns}. Their strengths lie in flexibility, simplicity, and the embedding of physical knowledge into the data-driven model, making them very popular for various applications \cite{faroughi2024physics}. 
However, due to the use of MLPs in their architectures, they also face several challenges such as high computational costs for multi-dimensional PDEs, large-scale spatiotemporal domain, sensitivity to hyperparameter tuning, and training difficulties such as vanishing gradients \cite{abbasi2024physics, cao2024surrogate, hu2024tackling, wang2020understanding}. In addition, resolving sharp transitions in highly dynamical systems often requires a lot of hand-tuning, which restricts their scalability \cite{costabal2024delta, liu2024discontinuity}.

To address the limitations of traditional PINNs based on MLPs, the scientific machine learning community has begun to explore replacements for MLPs, such as a prominent one known as the Kolmogorov-Arnold Network (KAN) \cite{ liu2024kan, KAN, shukla2024comprehensive}. KANs are inspired by the Kolmogorov-Arnold representation theorem that expresses multivariate functions as sums and compositions of univariate functions \cite{ braun2009constructive, kolmogorov1957representation}. Unlike MLPs, KANs replace fixed activation functions with learnable univariate functions, offering a more flexible and adaptive architecture.
Physics-informed Kolmogorov-Arnold Networks (PIKANs) \cite{shukla2024comprehensive} combine the strengths of KANs and PINNs, offering a major advance in scientific machine learning.
PIKANs combine the interpretability of KANs with the physics-guided learning capabilities of PINNs, making them highly effective in capturing the dynamics of systems governed by complex PDEs and boundary conditions \cite{howard2024finite, qiu2024relu, toscano2024kkans}.
Recent studies \cite{shukla2024comprehensive} have demonstrated that while vanilla KANs with B-splines struggle with robustness and efficiency, variants employing low-order orthogonal polynomials, such as radial basis functions (RBF) \cite{li2024kolmogorov, mostajeran2023radial}, wavelets \cite{bozorgasl2024wav, mostajeran2023novel}, or Chebyshev polynomials \cite{ss2024chebyshev} achieve competitive or superior performance than PINNs. This flexibility has led to the development of different implementations \cite{alsaadi2022control,  deepthi2023development, guo2024physics, mostajeran2024epi}. Despite these strengths, KANs are found to be computationally intensive for high-dimensional problems, and hence, they are slower compared to  MLPs \cite{shukla2024comprehensive}. Moreover, they are sensitive to the choice of basis functions and optimization strategies, which should be designed with care to achieve robust and efficient performance \cite{mostajeran2024epi}. These challenges point to the necessity for further advances in KAN-based approaches to make the most of their potential in scientific computing. Therefore, improving the performance and scalability of PIKANs is an active area of research to overcome their limitations and increase their applicability.
In this regard, \cite{guo2024physics} introduced ChebPIKANs, a model that includes Chebyshev polynomials as basis functions in their construction for enhanced approximation capabilities, stability, and extrapolation in the fluid dynamics problem, avoiding possible overfitting. Constructed upon KANs, ChebPIKAN retains all the benefits of locally learned representations based on the spline and overcomes generalization challenges in conventional neural networks. Moreover, the inclusion of physical information not only enhances the performance, but also provides interpretability to the model to some extent, giving it a more physically meaningful setting than the vanilla MLP architectures. Although ChebPIKANs attain better adaptability and robustness, the mathematical interpretation of their architecture remains a challenge due to the several algorithmic optimizations involved in the process.
\cite{wang2024kolmogorov} also tried to solve different forms of PDEs such as strong form, energy form, and inverse form with Kolmogorov–Arnold-informed neural network (KINN) which essentially involved replacing the MLP with KAN. Although KINN significantly outperforms MLP in terms of accuracy and convergence speed for numerous PDEs in computational solid mechanics, the authors discovered that because of the conflict between grid size and geometric complexity, KAN did not perform as well on complex geometric PDE problems.
More recently, \cite{jacob2024spikans} introduced Separable PIKANs (SPIKANs), which apply the principle of variable separation to reduce the computational complexity of KANs. By dedicating each dimension to an independent KAN, SPIKANs achieve remarkable training efficiency for high-dimensional problems without sacrificing accuracy. This decomposition distributes \(O(N)\) points to each KAN,  reducing memory usage and computation time. However, SPIKANs require a factorizable grid of collocation points rather than the unstructured point clouds typically used in PINNs, which can limit their applicability in some scenarios.
These investigations highlight the increasing initiatives to enhance PIKAN frameworks, tackling computational obstacles and broadening their applicability in various scientific domains.

Chebyshev polynomials are defined over \([-1, 1]\), which requires input normalization in cKAN layers to ensure stability. \cite{ss2024chebyshev} proposed a forward-pass normalization strategy for this purpose. Studies such as \cite{guo2024physics, hu2024tackling, ss2024chebyshev} emphasize that the absence of normalization often leads to training instability, including divergence and NaN loss after several iterations \cite{shukla2024comprehensive, ss2024chebyshev}. For example, \cite{shukla2024comprehensive} observed loss instability when solving high-frequency Helmholtz equations using the cPIKAN model.
In this work, we introduce Scaled-cPIKAN, a physics-informed neural network framework that synergizes Chebyshev polynomial-based representations with a scaling approach for spatial variables in governing equations. By combining the adaptive mathematical properties of KANs with the physics-driven principles of PINNs, Scaled-cPIKAN effectively captures high-frequency oscillations and intricate dynamics critical to multi-scale PDE modeling.
This architecture works great for complex PDE problems over extended spatial domains, where other approaches face some bottlenecks, such as poor convergence and loss of accuracy. 
In particular, our approach highlights the critical role of spatial domain scaling in addressing the computational inefficiencies and accuracy limitations associated with solving PDEs over extended spatial domains with highly oscillatory dynamics. Such a scaling strategy remarkably improves the accuracy of the approximation and increases the rate of loss convergence, as shown in four different test cases.
The subsequent sections of this paper are structured as follows. Section \ref{Sec.Pre} summarizes foundational concepts such as MLPs, cKANs, and PINNs. In Section \ref{Sec.method}, we introduce Scaled-cPIKANs, and demonstrate its efficacy against other methods in Section \ref{Sec.res} using four benchmark problems. 
Finally, the main conclusions are drawn in Section \ref{Sec.Conclusion}.

%%%=============================================================
\section{Preliminaries}\label{Sec.Pre}

Partial differential equations (PDEs) are fundamental mathematical tools that model most physical, chemical, and biological phenomena. They express the dynamics of complex systems across spatiotemporal domains in terms of relations between spatial and temporal derivatives. They are generally bounded by initial and boundary conditions to develop a well-posed formulation of the problem \cite{evans2022partial, hadamard1902problemes}.
A general PDE defined over a spatial domain \(\Omega \subset \mathbb{R}^d\) and a temporal interval $t \in [0, T]$  can be stated as follows,
\begin{equation}\label{Eq.GeneralPDE}
    \begin{array}{ll}
         \mathcal{N}[u(\boldsymbol{x}, t); \Lambda] = f(\boldsymbol{x}, t), &  \boldsymbol{x}\in \Omega,\: t\in (0,T),\\[4mm]
        u(\boldsymbol{x}, 0) = u_0(\boldsymbol{x}), & \boldsymbol{x}\in \Omega,\\[4mm]
        \mathcal{B}[u(\boldsymbol{x}, t)] = u_b(\boldsymbol{x}, t), & \boldsymbol{x}\in \partial\Omega,\: t\in (0,T),
    \end{array}
\end{equation}
where 
$\boldsymbol{x} = (x_1, x_2, \ldots, x_d)$ represents the spatial coordinates and $t$ is the temporal variable,
$u$ is the unknown attribute,
$\mathcal{N}[\cdot \: ; \Lambda]$ is a general spatial-temporal differential operator parametrized by $\Lambda$, which may be nonlinear,
$\mathcal{B}$ is a boundary operator,
and
$f$ is a source term defined over $\Omega \times [0,T]$.
In Eq.~\eqref{Eq.GeneralPDE},  $u_0$ specifies the initial condition at $t=0$ and $u_b$ specifies the boundary condition at the spatial boundary $\partial\Omega$.
This formulation provides a general framework for analyzing and solving PDEs across various applications. 
The structure of Eq.~\eqref{Eq.GeneralPDE} ensures flexibility in incorporating complex dynamics, making it suitable for modeling nonlinear, multiscale, or coupled systems. 
In the forward PDE problems, the objective is to determine the solution \(u\) given the governing equations, parameters, and initial/boundary conditions, where the parameter \(\Lambda\) is also known. In contrast, inverse PDE problems focus on identifying unknown quantities, such as the parameter \(\Lambda\) or source terms \(f\), from observed data (i.e., a partial solution provided for \(u\)). Inverse problems are often ill-posed \cite{tikhonov1977solutions}, necessitating additional regularization techniques or constraints to obtain meaningful and stable solutions. In neural-driven PDE-solvers, various neural network architectures are employed to approximate these unknowns, a few of which, used in this work, are briefly reviewed below. 
%%%=============================================================
\subsection{Multilayer Perceptron}\label{Sec.MLP}

A Multilayer Perceptron (MLP) is a key building block of feed-forward artificial neural networks \cite{bengio2017deep, bishop2006pattern, sazli2006brief, svozil1997introduction}. It consists of multiple layers: an input layer, one or several hidden layers, and an output layer. The main role of an MLP is to approximate complex, nonlinear relationships in data by a series of transformations, in which activation functions play an important role.
The analysis performed by an MLP can be viewed as a series of nested mappings, performed layer after layer. Denote \(\boldsymbol{\xi}\) as the input vector and \(L\) as the total number of layers,  then, the output \( y \) of the MLP can be expressed as
\begin{equation}
    y_{_\text{MLP}}(\boldsymbol{\xi}) = \mathcal{F}^{(L)} \circ \mathcal{F}^{(L-1)} \circ \dots \circ \mathcal{F}^{(1)}(\boldsymbol{\xi}),
\end{equation}
where each \(\mathcal{F}^{(l)}\), corresponding to the \(l\)-th layer, represents the transformation involving the weights, biases, and activation function of that layer. The operation in the \(l\)-th layer can be expressed as,
\begin{equation}\label{Eq.Layer}
    \mathcal{F}^{(l)}(\boldsymbol{\xi}^{(l)}) = \sigma\left( \boldsymbol{W}^{(l)} \boldsymbol{\xi}^{(l-1)} + \boldsymbol{b}^{(l)} \right),
\end{equation}  
where $\boldsymbol{\xi}^{(l)}$ represents the input vector to the \(l\)-th layer, with $\boldsymbol{\xi}^{(0)} = \boldsymbol{\xi}$ as the initial input to the network. 
In Eq.~\eqref{Eq.Layer}, \(\boldsymbol{W}^{(l)}\) and \(\boldsymbol{b}^{(l)}\) refer to the weights and biases of the \(l\)-th layer, respectively, and \(\sigma\) is the activation function that can take different fixed forms, e.g., sin, tanh,  \cite{faroughi2023physics}.   
The set of trainable parameters in an MLP can be denoted by  
$
\boldsymbol{\theta}_{_\text{MLP}} = \{ \boldsymbol{W}^{(l)}, \boldsymbol{b}^{(l)} \}.
$  
The total number of parameters in the network is approximately proportional to \(O(N_l N_n^2)\), where \(N_l\: (= L-1)\) denotes the number of hidden layers, and \(N_n\) is the number of neurons per hidden layer \cite{shukla2024comprehensive}. 
The multilayer architecture of MLPs allows subsequent layers to take the output of other layers as input, to model very complicated relationships.
The Universal Approximation Theorem \cite{cybenko1989approximation} guarantees that if there are enough neurons and with appropriate activation functions, MLPs can approximate virtually any continuous function. This capability makes MLPs powerful tools in machine learning.
%%%=============================================================
\subsection{Chebyshev-based Kolmogorov-Arnold Networks}\label{Sec.KAN}

Kolmogorov-Arnold networks (KANs) are another set of architecture developed based on the Kolmogorov-Arnold representation theorem \cite{liu2024kan, KAN}. This theorem guarantees that any continuous multivariate function on a bounded domain can be represented by a composition of continuous univariate functions and additive operations. 
Formally, for \( h: [0,1]^d \to \mathbb{R} \),
\begin{equation}\label{Eq.SimpleKAN}
h(\boldsymbol{\xi}) = \sum_{j=0}^{2d} \Phi_j\left(\sum_{i=1}^d \phi_{i,j}(\xi_i)\right),
\end{equation}
where \(\phi_{i,j}: [0,1] \to \mathbb{R}\) and \(\Phi_j: \mathbb{R} \to \mathbb{R}\) are univariate functions. Variations in the forms of these functions define different KAN architectures. KANs generalize this representation by organizing computations across \(L\) layers, where each layer consists of univariate transformations and linear combinations,
\begin{equation}\label{Eq.yKAN}
y_\text{KAN}(\boldsymbol{\xi}) = (\boldsymbol{\Phi}_L \circ \cdots \circ \boldsymbol{\Phi}_1)(\boldsymbol{\xi}),
\end{equation}
where  \(\boldsymbol{\Phi}_l\) is a KAN layer defined as a matrix of 1-dimensional functions,
\begin{equation}\label{Eq.KANlayer}
    \boldsymbol{\Phi}_l = \{\phi_{l, i,j}\},\:\: l=1, \cdots , L,\:\: i=1,\cdots , n_{l}, \:\: j=1,\cdots , n_{l+1}. 
\end{equation}

In Eq.~\eqref{Eq.KANlayer}, \(n_{l}\) and \(n_{l+1}\) are the inputs and outputs dimension, respectively.
% with each layer \(\boldsymbol{\Phi}_l\) represented as a matrix of univariate functions \(\phi_{l,i,j}\).
KANs often employ basis functions in the form of splines to enhance flexibility,
\begin{equation}
\phi(\xi) = w_b b(\xi) + w_s \text{spline}(\xi), \quad b(\xi) = \frac{\xi}{1 + e^{-\xi}}, \quad \text{spline}(\xi) = \sum_i c_i B_i(\xi),
\end{equation}
where \(B_i(\xi)\) are spline basis functions defined by the grid size \(g\) and polynomial degree \(k\). However, increasing \(g\) greatly increases the number
of trainable parameters $|\boldsymbol{\theta}|_\text{KAN} \sim O(N_l N_n^2 (k + g))$,
where \(N_l\) is the number of layers, \(N_n\) the number of neurons per layer, \(g\) the grid size, and \(k\) the polynomial degree. One way to address this computational issue is to employ Chebyshev polynomials, a fundamental tool in approximation theory and numerical analysis,  as an orthogonal basis  \cite{rivlin2020chebyshev, schmidt2021kolmogorov}.  
Chebyshev polynomials are defined recursively as,
\begin{equation}
\begin{aligned}
T_0(\xi) &= 1, \quad T_1(\xi) = \xi, \\
T_n(\xi) &= 2\xi T_{n-1}(\xi) - T_{n-2}(\xi), \quad n \geq 2.
\end{aligned}
\end{equation}
following which, the univariate functions in Chebyshev-based KANs (cKANs) are represented as,
\begin{equation}
\phi(\xi) = \sum_i c_i T_i(\xi),
\end{equation}
where \(c_i\) are the trainable coefficients. This reduces the number of trainable parameters to 
$|\boldsymbol{\theta}|_\text{cKAN} \sim O(N_l N_n^2 k)$.
To improve stability while training, cKANs restrict their inputs in the range \([-1, 1]\). That is accomplished by applying the  \(\tanh\)  after each layer \cite{guo2024physics, hu2024tackling, ss2024chebyshev}, changing the corresponding forward pass, described by Eq.~\eqref{Eq.yKAN} as,
\begin{equation}  
y_\text{cKAN}(\boldsymbol{\xi}) = (\boldsymbol{\Phi}_L \circ \tanh \circ \cdots \circ \tanh \circ ~ \boldsymbol{\Phi}_1 \circ \tanh)(\boldsymbol{\xi}).  
\end{equation}

%%%=============================================================
\subsection{Physics-informed  Networks}\label{Sec.PI}

The concept behind physics-informed networks for solving the general form of PDE presented in Eq.~\eqref{Eq.GeneralPDE} aims to utilize the neural networks' expressive capabilities to approximate the solution, e.g.,  PINNs using a MLP architecture denoted as \(u(\boldsymbol{\xi}; \boldsymbol{\theta}_{\text{MLP}})\) or PIKANs and cPIKANs using a KAN architecture denoted as \(u(\boldsymbol{\xi}; \boldsymbol{\theta}_{\text{KAN}})\). 
Physics-informed networks approximate solutions by directly integrating governing equations, boundary conditions, and initial conditions into the loss function, thereby utilizing the underlying physics. In this way, it is guaranteed that the approximated solution fits the observed data, if available,  and satisfies the physical laws dictated by the governing equations. This approach introduces an additional layer for regularization, hence making the technique highly feasible in solving both the inverse and forward problems.

In physics-informed networks, the training process utilizes two primary types of data. The first type consists of  collocation points, denoted as \(\mathcal{T}_{\text{res}} = \{\boldsymbol{\xi}^{\text{res}}_i \in \Omega \times (0,T)\}_{i=1}^{N_{\text{res}}}\), which lie in the spatial domain \(\Omega\) and the time interval \((0,T)\). The second type consists of initial, boundary, and measurement data.
The  initial data is represented as \(\mathcal{T}_{\text{init}} = \{\boldsymbol{\xi}_j^{\text{init}}, u_0(\boldsymbol{\xi}_j^{\text{init}})\}_{j=1}^{N_{\text{init}}}\), where the input points \(\boldsymbol{\xi}_j^{\text{init}} = (\boldsymbol{x}_j, 0)\) correspond to the initial time (\(t=0\)) with \(\boldsymbol{x}_j \in \Omega\) and \(u_0(\boldsymbol{\xi}_j^{\text{init}})\) representing the prescribed initial state.
The  boundary data  is represented as \(\mathcal{T}_{\text{bc}} = \{\boldsymbol{\xi}_j^{\text{bc}}, u_b(\boldsymbol{\xi}_j^{\text{bc}})\}_{j=1}^{N_{\text{bc}}}\), where \(\boldsymbol{\xi}_j^{\text{bc}} = (\boldsymbol{x}_j, t_j) \in \partial \Omega \times (0,T)\) and \(u_b(\boldsymbol{\xi}_j^{\text{bc}})\) denote the prescribed boundary values. 
Finally, the measurement data, used mainly in solving inverse problems,  is captured by \(\mathcal{T}_{\text{meas}} = \{\boldsymbol{\xi}_j^{\text{meas}}, u_{\text{meas}}(\boldsymbol{\xi}_j^{\text{meas}})\}_{j=1}^{N_{\text{meas}}}\), where  \(\boldsymbol{\xi}_j^{\text{meas}} = (\boldsymbol{x}_j, t_j) \in \Omega \times (0,T)\) 
and \(u_{\text{meas}}(\boldsymbol{\xi}_j^{\text{meas}})\) correspond to the measurement locations and values in the computational domain. These  datasets are used to train networks by minimizing a loss function that incorporates both the residuals of the governing equations and the error between the predicted solution and the given data. The general form of the loss function is thus defined as, 
\begin{equation}\label{Eq.Loss}
\mathcal{L}(\boldsymbol{\theta}) = \lambda_{\text{res}} \mathcal{L}_{\text{res}}(\mathcal{T}_{\text{res}}; \boldsymbol{\theta}) + \lambda_{\text{data}} \mathcal{L}_{\text{data}}(\mathcal{T}_{\text{data}}; \boldsymbol{\theta}),
\end{equation}
where \(\mathcal{T}_{\text{data}} = \mathcal{T}_{\text{init}} \cup \mathcal{T}_{\text{bc}} \cup \mathcal{T}_{\text{meas}}\), and \(\lambda_{\text{res}}\) and \(\lambda_{\text{data}}\) are scalar coefficients that control the relative importance of the physics residual loss \(\mathcal{L}_{\text{res}}\) and the data fitting loss \(\mathcal{L}_{\text{data}}\).

%%%=============================================================
\section{Domain Scaling in cPIKANs (scaled-cPIKANs)}\label{Sec.method}

When dealing with extended spatial computation domains and problems containing high-frequency components, normalizing the inputs to \([-1, 1]\) in each cKAN layer can make it more stable but may not avoid instabilities. The Helmholtz equation for high wave numbers is one problem that has faced loss divergence in cPIKAN models after a few iterations \cite{hu2024tackling}.
The main idea of our proposed method is to scale spatial variables in PDEs to the standardized interval \([-1, 1]\), as required by Chebyshev polynomials. This approach enables solving the scaled PDE, which has been reformulated to operate over a scaled spatial domain, instead of the original one presented in Eq.~\eqref{Eq.GeneralPDE}. This transformation not only ensures compatibility with Chebyshev-based methods, but also facilitates numerical stability and efficiency.
Consider the scaled PDE defined as follows,
\begin{equation}\label{Eq.ScaledPDE}
    \begin{array}{ll}
        \mathcal{N}^{\text{s}}[\tilde{u}(\tilde{\boldsymbol{x}}, t); \Lambda] = \tilde{f}(\tilde{\boldsymbol{x}}, t), &  \tilde{\boldsymbol{x}}\in \tilde{\Omega},\: t\in (0,T),\\[4mm]
        \tilde{u}(\tilde{\boldsymbol{x}}, 0) = \tilde{u}_0(\tilde{\boldsymbol{x}}), & \tilde{\boldsymbol{x}}\in \tilde{\Omega},\\[4mm]
        \mathcal{B}^{s}[\tilde{u}(\tilde{\boldsymbol{x}}, t)] = \tilde{u}_b(\tilde{\boldsymbol{x}}, t), & \tilde{\boldsymbol{x}}\in \partial\tilde{\Omega},\: t\in (0,T),
    \end{array}
\end{equation}
where \(\tilde{\Omega} = [-1, 1]^d\), and the tilde over the functions indicates that these functions are defined on the scaled domain, distinguishing them from their counterparts in the original domain. The operators \(\mathcal{N}^{\text{s}}\) and \(\mathcal{B}^{\text{s}}\) represent the scaled versions of the original differential and boundary operators, respectively, with the superscript \(s\) signifying the scaling transformation.

In the scaled domain, the training data must be transformed accordingly to maintain consistency with the scaled PDE formulation. The set of residual collocation points in the scaled domain is denoted by \(\mathcal{T}_{\text{res}}^{\text{s}} = \{\tilde{\boldsymbol{\xi}}^{\text{res}}_i \in \tilde{\Omega} \times (0, T)\}_{i=1}^{N_{\text{res}}}\), where \(\tilde{\boldsymbol{\xi}}^{\text{res}}_i = (\tilde{\boldsymbol{x}}_i, t_i)\). Similarly, the initial condition data points are represented as \(\mathcal{T}_{\text{init}}^{\text{s}} = \{\tilde{\boldsymbol{\xi}}_j^{\text{init}}, \tilde{u}_0(\tilde{\boldsymbol{\xi}}_j^{\text{init}})\}_{j=1}^{N_{\text{init}}}\), where \(\tilde{\boldsymbol{\xi}}_j^{\text{init}} = (\tilde{\boldsymbol{x}}_j, 0)\), and \(\tilde{u}_0(\tilde{\boldsymbol{\xi}}_j^{\text{init}})\) correspond to the prescribed initial state in the scaled domain. For boundary conditions, the data set is expressed as \(\mathcal{T}_{\text{bc}}^{\text{s}} = \{\tilde{\boldsymbol{\xi}}_j^{\text{bc}}, \tilde{u}_b(\tilde{\boldsymbol{\xi}}_j^{\text{bc}})\}_{j=1}^{N_{\text{bc}}}\), where \(\tilde{\boldsymbol{\xi}}_j^{\text{bc}} = (\tilde{\boldsymbol{x}}_j, t_j) \in \partial\tilde{\Omega} \times (0, T)\), and \(\tilde{u}_b(\tilde{\boldsymbol{\xi}}_j^{\text{bc}})\) correspond to the boundary values in the scaled spatial domain. In addition, the measurement data are reformulated as \(\mathcal{T}_{\text{meas}}^{\text{s}} = \{\tilde{\boldsymbol{\xi}}_j^{\text{meas}}, \tilde{u}_{\text{meas}}(\tilde{\boldsymbol{\xi}}_j^{\text{meas}})\}_{j=1}^{N_{\text{meas}}}\), where \(\tilde{\boldsymbol{\xi}}_j^{\text{meas}} = (\tilde{\boldsymbol{x}}_j, t_j) \in \tilde{\Omega} \times (0, T)\), and \(\tilde{u}_{\text{meas}}(\tilde{\boldsymbol{\xi}}_j^{\text{meas}})\) represent the scaled locations and corresponding observed values. Accordingly, the total training data in the scaled domain is defined as \(\mathcal{T}_{\text{data}}^{\text{s}} = \mathcal{T}_{\text{init}}^{\text{s}} \cup \mathcal{T}_{\text{bc}}^{\text{s}} \cup \mathcal{T}_{\text{meas}}^{\text{s}}\). 
To address the challenges posed by spatial variable scaling in PDEs, the loss function is reformulated to reflect the scaled spatial domain while maintaining consistency with the governing physical laws. The scaled loss function is expressed as,  
\begin{equation}\label{Eq.SLoss}
    \mathcal{L}^{\text{s}}(\boldsymbol{\theta}) = \lambda_{\text{res}} \mathcal{L}_{\text{res}}^{\text{s}}(\boldsymbol{\theta}) + \lambda_{\text{data}} \mathcal{L}_{\text{data}}^{\text{s}}(\boldsymbol{\theta}),
\end{equation}
where \(\mathcal{L}_{\text{res}}^{\text{s}}\) represents the residual loss on the scaled domain, and \(\mathcal{L}_{\text{data}}^{\text{s}}\) incorporates the scaled contributions of initial conditions, boundary conditions, and measurement data.  
The residual loss in the scaled domain is defined as,
\begin{equation}\label{Eq.SLossRes}
    \mathcal{L}_{\text{res}}^{\text{s}}(\mathcal{T}_{\text{res}}^{\text{s}}; \boldsymbol{\theta}) = \frac{1}{N_{\text{res}}} \sum_{i=1}^{N_{\text{res}}} \vert\mathcal{N}^{\text{s}}[\tilde{u}(\tilde{\boldsymbol{\xi}}_i^{\text{res}};\Lambda, \boldsymbol{\theta})] - \tilde{f}(\tilde{\boldsymbol{\xi}}_i^{\text{res}})\vert^2,
\end{equation}  
where \(\mathcal{N}^{\text{s}}\) is the differential operator adjusted for the scaled domain, and \(\tilde{f}\) is the source term in the scaled domain.   
The data loss, \(\mathcal{L}_{\text{data}}(\boldsymbol{\theta})\), in Eq.~\eqref{Eq.SLoss}, can be decomposed into three distinct components, each addressing different aspects of the problem: initial conditions, boundary conditions, and measurement data. 
Formally, the scaled data loss is defined as,
\begin{equation}\label{Eq.SLossData}
    \mathcal{L}_{\text{data}}^{\text{s}}(\boldsymbol{\theta}) = 
\lambda_{\text{init}} \mathcal{L}_{\text{init}}^{\text{s}}(\boldsymbol{\theta}) + 
\lambda_{\text{bc}} \mathcal{L}_{\text{bc}}^{\text{s}}(\boldsymbol{\theta}) +
\lambda_{\text{meas}} \mathcal{L}_{\text{meas}}^{\text{s}}(\boldsymbol{\theta}),
\end{equation} 
where \(\mathcal{L}^{\text{s}}_{\text{init}}(\boldsymbol{\theta})\) penalizes deviations from the specified initial conditions, 
\(\mathcal{L}^{\text{s}}_{\text{bc}}(\boldsymbol{\theta})\) enforces the satisfaction of boundary conditions, and
\(\mathcal{L}^{\text{s}}_{\text{meas}}(\boldsymbol{\theta})\) ensures consistency with any additional measurement data provided.
\(\lambda_{\text{init}}\) and \(\lambda_{\text{bc}}\) are set to zero if initial or boundary conditions are not relevant or available for the specific problem.  Moreover, \(\lambda_{\text{meas}}\) is typically set to one for inverse problems, where measurement data guide the solution, and to zero for forward problems, where the focus lies solely on satisfying the governing equations and initial/boundary conditions.
In Eq.~\eqref{Eq.SLossData}, each loss term is defined as follows,
\begin{equation}\label{Eq.SLossINIT}
    \mathcal{L}^{\text{s}}_{\text{init}}(\mathcal{T}^{\text{s}}_{\text{init}}; \boldsymbol{\theta}) = \frac{1}{N_{\text{init}}} \sum_{j=1}^{N_{\text{init}}} \vert\tilde{u}(\tilde{\boldsymbol{\xi}}_j^{\text{init}}; \boldsymbol{\theta}) - \tilde{u}_0(\tilde{\boldsymbol{\xi}}_j^{\text{init}})\vert^2,
\end{equation}
\begin{equation}\label{Eq.SLossBC}
    \mathcal{L}^{\text{s}}_{\text{bc}}(\mathcal{T}^{\text{s}}_{\text{bc}}; \boldsymbol{\theta}) = \frac{1}{N_{\text{bc}}} \sum_{j=1}^{N_{\text{bc}}} \vert\mathcal{B}^{\text{s}}[\tilde{u}(\tilde{\boldsymbol{\xi}}_j^{\text{bc}}; \boldsymbol{\theta})] - \tilde{u}_b(\tilde{\boldsymbol{\xi}}_j^{\text{bc}}) \vert^2,
\end{equation}
and
\begin{equation}\label{Eq.SLossMEAS}
    \mathcal{L}^{\text{s}}_{\text{meas}}(\mathcal{T}^{\text{s}}_{\text{meas}};\boldsymbol{\theta}) = \frac{1}{N_{\text{meas}}} \sum_{j=1}^{N_{\text{meas}}} \vert \tilde{u}(\tilde{\boldsymbol{\xi}}_j^{\text{meas}}; \boldsymbol{\theta}) - \tilde{u}_{\text{meas}}(\tilde{\boldsymbol{\xi}}_j^{\text{meas}})\vert^2.
\end{equation}

In forward problems, the objective is to find the optimal network parameters \(\boldsymbol{\theta}^*\) such that the predicted solution \(\tilde{u}(\tilde{\boldsymbol{\xi}}; \boldsymbol{\theta}^*)\) satisfies the PDE, initial conditions, and boundary conditions. This optimization can be expressed as
\begin{equation}
    \boldsymbol{\theta}^* = \arg\min_{\boldsymbol{\theta}} \mathcal{L}^{\text{s}}(\boldsymbol{\theta}).
\end{equation}

For inverse problems, the aim is to infer unknown PDE parameters, source terms, or boundary conditions from measurement data. 
The loss function for such problems is tailored to emphasize measurement data, with \(\lambda_{\text{meas}} = 1\) in Eq.~\eqref{Eq.SLossData}, reflecting the increased reliance on observational information. 
The optimization problem for inverse problems is formulated as,
\begin{equation}
    \boldsymbol{\theta}^*, \Lambda^* = \arg\min_{\boldsymbol{\theta}, \Lambda} \mathcal{L}^{\text{s}}(\boldsymbol{\theta}, \Lambda),
\end{equation}
where \(\Lambda\) represents the set of unknown PDE parameters.
Through this process, the network simultaneously learns the optimal network parameters 
\(\boldsymbol{\theta}^*\) and estimates the unknown PDE parameters 
\(\Lambda^*\). The training process employs gradient-based optimization methods, with gradients of $\mathcal{L}^{\text{s}}(\boldsymbol{\theta})$ computed using automatic differentiation. This enables efficient and accurate parameter updates, converging to a solution \(\tilde{u} \approx \tilde{u}(\boldsymbol{\xi}; \boldsymbol{\theta}^*)\).

%%%=============================================================
\section{Results and Discussion}\label{Sec.res}
In this section, we perform several numerical experiments to test the efficiency of the proposed method, Scaled-cPIKANs, compared to other  methods.
The numerical experiments include the diffusion equation across multiple domain sizes, the Helmholtz equation under varying PDE parameter conditions, the Allen-Cahn equation in both supervised and unsupervised contexts, as well as the reaction-diffusion equation in forward and inverse problem settings.
To distinguish between methods in comparisons, we refer to the spatial domain-scaled versions of cPIKAN and PINN as ``Scaled-cPIKAN'' and ``Scaled-PINN'', respectively, while the unscaled (i.e., vanilla) versions are simply referred to as ``cPIKAN'' and ``PINN''. Our code is implemented in PyTorch on a machine equipped with an NVIDIA RTX 6000 Ada GPU. All neural networks in this work are trained using the Adam optimizer with an initial learning rate of \(10^{-3}\) unless otherwise specified. The total number of training epochs is indicated in each subsection. 
In our comparison, to quantify the accuracy of the schemes' predictions against the ground truth solutions, we report the relative \(\mathcal{L}^2\) error, defined as,
\begin{equation}
\label{Eq.RelL2Error}
\mathrm{RE} 
= \frac{
    \sqrt{
      \displaystyle \sum_{\boldsymbol{\xi} \in \mathcal{T}^*}
        \bigl|\,u_{\mathrm{E}}(\boldsymbol{\xi}) - u_{\mathrm{P}}(\boldsymbol{\xi})\bigr|^{2}
    }
  }{
    \sqrt{
      \displaystyle \sum_{\boldsymbol{\xi} \in \mathcal{T}^*}
        \bigl|\,u_{\mathrm{E}}(\boldsymbol{\xi})\bigr|^{2}
    }
  },
\end{equation}
where \(u_{\text{E}}\) is the ground truth solution, and \(u_{\text{P}}\) is the approximated solution obtained from one of the four methods: Scaled-cPIKAN, Scaled-PINN, cPIKAN, and PINN. The test set \(\mathcal{T}^*\) consists of collocation points in the spatiotemporal domain that are different from those used in the training set. 

%%%=============================================================
%%%=============================================================
\subsection{Diffusion Equation }\label{Exam.DiffEqu}
%\begin{example}\label{Exam.DiffEqu}
As the first experiment, we aim to solve the one-dimensional diffusion equation with Dirichlet boundary conditions and an initial condition given as,
\begin{equation}
    % u_{t} - D \, u_{xx} = 0, \quad x \in [-M, M], \: t \in (0, T], 
    \begin{array}{ll}
        u_{t}(x,t) - D\: u_{xx}(x,t) = 0, &  x \in [-M, M], t\in (0, T],\\[3mm]
        u(-M,t) = u(M,t) = 0, &t\in (0, T] ,\\[3mm]
        u(x, 0) = \sin(\pi x), &x \in [-M, M],
\end{array}
\end{equation}
 where $D$ is the diffusion coefficient, M indicates the half-width of the spatial domain, and T represents the final time. 
The ground truth solution for this problem is $u(x, t) = \sin(\pi x) \exp(-\pi^2 D t)$. We solve this example with four methods, Scaled-cPIKAN, Scaled-PINN, 
cPIKAN, and PINN for spatial domain sizes of \(M = 2, 4, 6\) and the final time $T=1$. The residual and data weights are set to \((\lambda_{\text{res}}, \lambda_{\text{data}}) = (0.01, 1.0)\), and the number of residual and data points are \((N_{\text{res}}, N_{\text{data}}) = (2000, 800)\). 
Since this problem is solved in a forward setting with no unknown PDE parameters, the data points consist only of the initial conditions and randomly selected boundary points at various time instances.

For each model, the number of hidden layers, $N_l$, the number of neurons per layer, $N_n$, and the degree of Chebyshe polynomial, $k$, are selected such that they contain almost the same amount of training parameters $\vert \boldsymbol{\theta}\vert$. In Table \ref{Tab.HeatSetting}, the configurations for each model and the associated \(\mathcal{L}^2\) errors for solving the diffusion equation with a diffusion coefficient of \(D = 0.1\) are presented.
From Table \ref{Tab.HeatSetting}, one can observe that Scaled-cPIKAN has the smallest relative errors on different sizes of the domain. The Scaled-PINN is more accurate than the vanilla PINN but worse than the Scaled-cPIKAN. Moreover, it is important to note that the efficacy of both PINN and cPIKAN decreases as the value of \(M\) increases. For example, at $M = 6$, the prediction given by Scaled-cPIKAN yields a relative \(\mathcal{L}^2\) error that is lower by 91.3\%, 96.2\%, and 97.5\% when compared to Scaled-PINN, PINN, and cPIKAN, respectively. This comparison illustrates the effectiveness of scaled methods, particularly Scaled-cPIKAN for extended spatial domains, when solving the diffusion equation.
%%%=============================================================
\begin{table}[!h]
\centering
\caption{\label{Tab.HeatSetting} 
Network configurations, number of parameters, and associated relative \(\mathcal{L}^2\) errors for solving the diffusion equation (Example \ref{Exam.DiffEqu}) for spatial domain sizes of \(M = 2, 4, 6\) and the final time $T=1$ assuming \(D = 0.1\), where \((\lambda_{\text{res}}, \lambda_{\text{data}}) = (0.01, 1.0)\) and \((N_{\text{res}}, N_{\text{data}}) = (2000, 800)\).}
\renewcommand{\arraystretch}{1.2}
\setlength{\tabcolsep}{12pt}
\begin{tabular}{c|cc|ccc}
\toprule
% \multicolumn{6}{l}{} \\
% \multicolumn{6}{l}{} \\
% \midrule
 &  &  & \multicolumn{3}{c}{Relative $\mathcal{L}^2$ Error} \\
Method & $(N_l, N_n, k)$ & $\vert \boldsymbol{\theta}\vert$ & $M=2$ & $M=4$ & $M=6$\\
\midrule
Scaled-cPIKAN & (2, 8, 5) & 462 &$\mathbf{4.3 \times 10^{-3}}$ & $\mathbf{1.3 \times 10^{-2}}$ & $\mathbf{2.1 \times 10^{-2}}$ \\
Scaled-PINN & (2, 19, -) & 456 & $1.2 \times 10^{-2}$& $3.1\times 10^{-2}$& $2.4\times 10^{-1}$\\
cPIKAN & (2, 8, 5) & 462 &$1.3\times 10^{-2}$ & $1.9\times 10^{-1}$& $8.4\times 10^{-1}$\\
PINN & (2, 19, -) & 456 & $2.0\times 10^{-2}$&$9.5\times 10^{-2}$&$5.5\times 10^{-1}$\\
%\midrule
\bottomrule
\end{tabular}
\end{table}
% % % ===========================================
\begin{figure}[!h]
    \centering
    \includegraphics[width=0.95\linewidth]{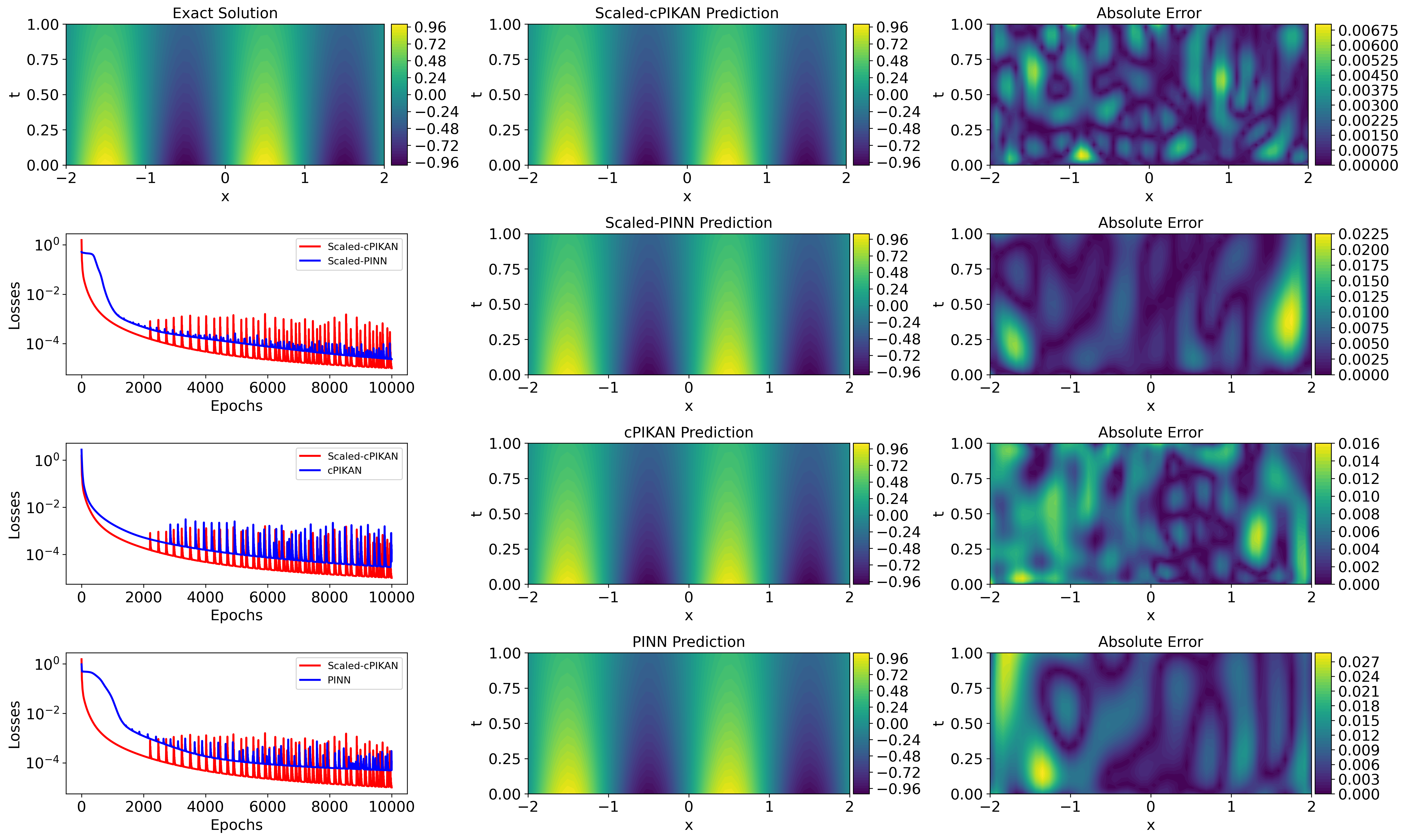}
    \caption{Comparison of the solutions predicted for the diffusion equation (Example \ref{Exam.DiffEqu}) in \([-2, 2] \times [0, 1]\). Each figure displays the outcomes of Scaled-cPIKAN, Scaled-PINN, cPIKAN, and PINN (from top to bottom). From left to right in each row: the loss function, the predicted solution \(u\), and the absolute error of the prediction are shown. The first plot in the top row shows the ground truth solution.}
    \label{fig:Heat2}
\end{figure}

\begin{figure}[!h]
    \centering
    \includegraphics[width=0.95\linewidth]{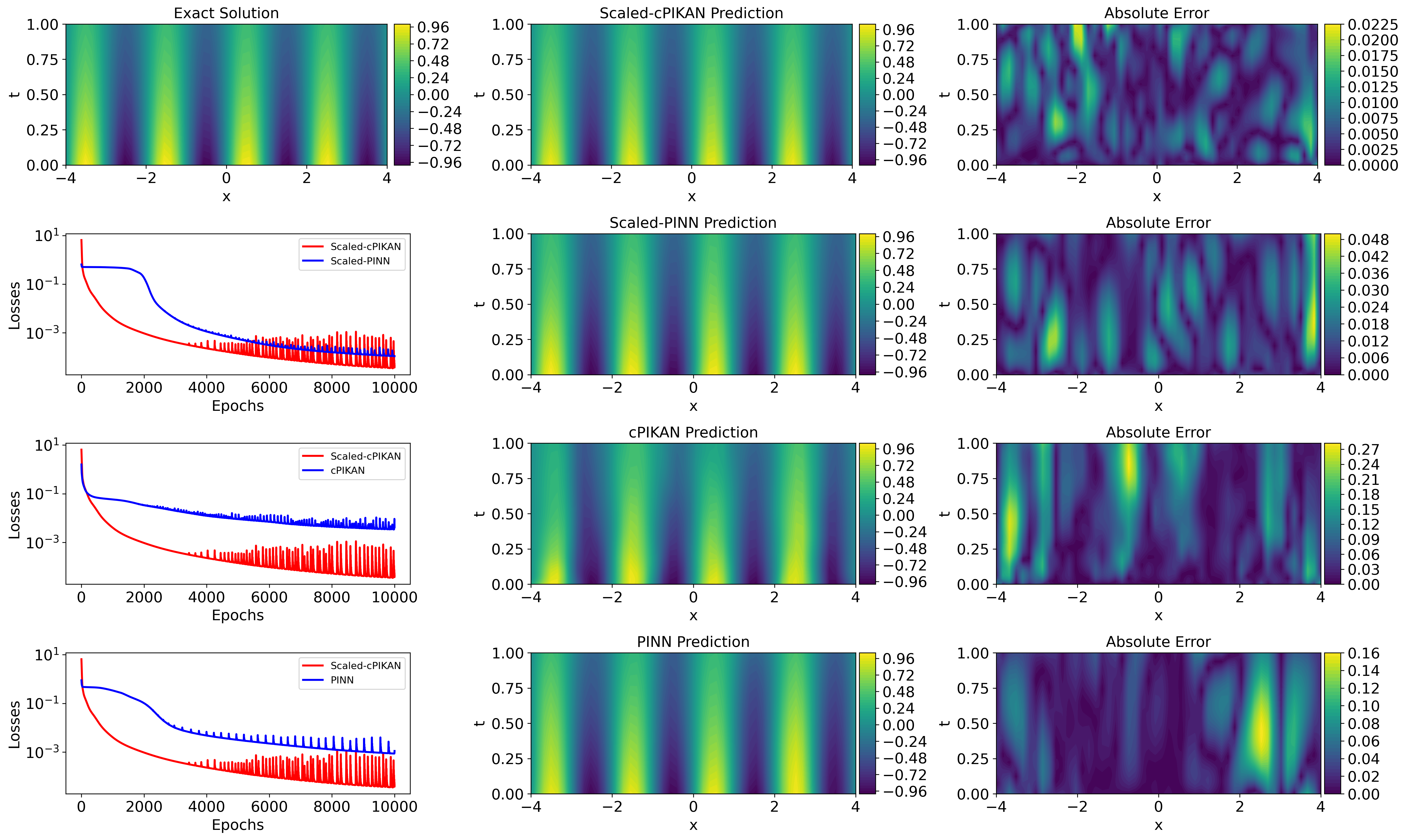}
    \caption{Comparison of the solutions predicted for the diffusion equation (Example \ref{Exam.DiffEqu}) in \([-4, 4] \times [0, 1]\). Each figure displays the outcomes of Scaled-cPIKAN, Scaled-PINN, cPIKAN, and PINN (from top to bottom). From left to right in each row: the loss function, the predicted solution \(u\), and the absolute error of the prediction are shown. The first plot in the top row shows the ground truth solution.}
    \label{fig:Heat4}
\end{figure}

\begin{figure}[!h]
    \centering
    \includegraphics[width=0.95\linewidth]{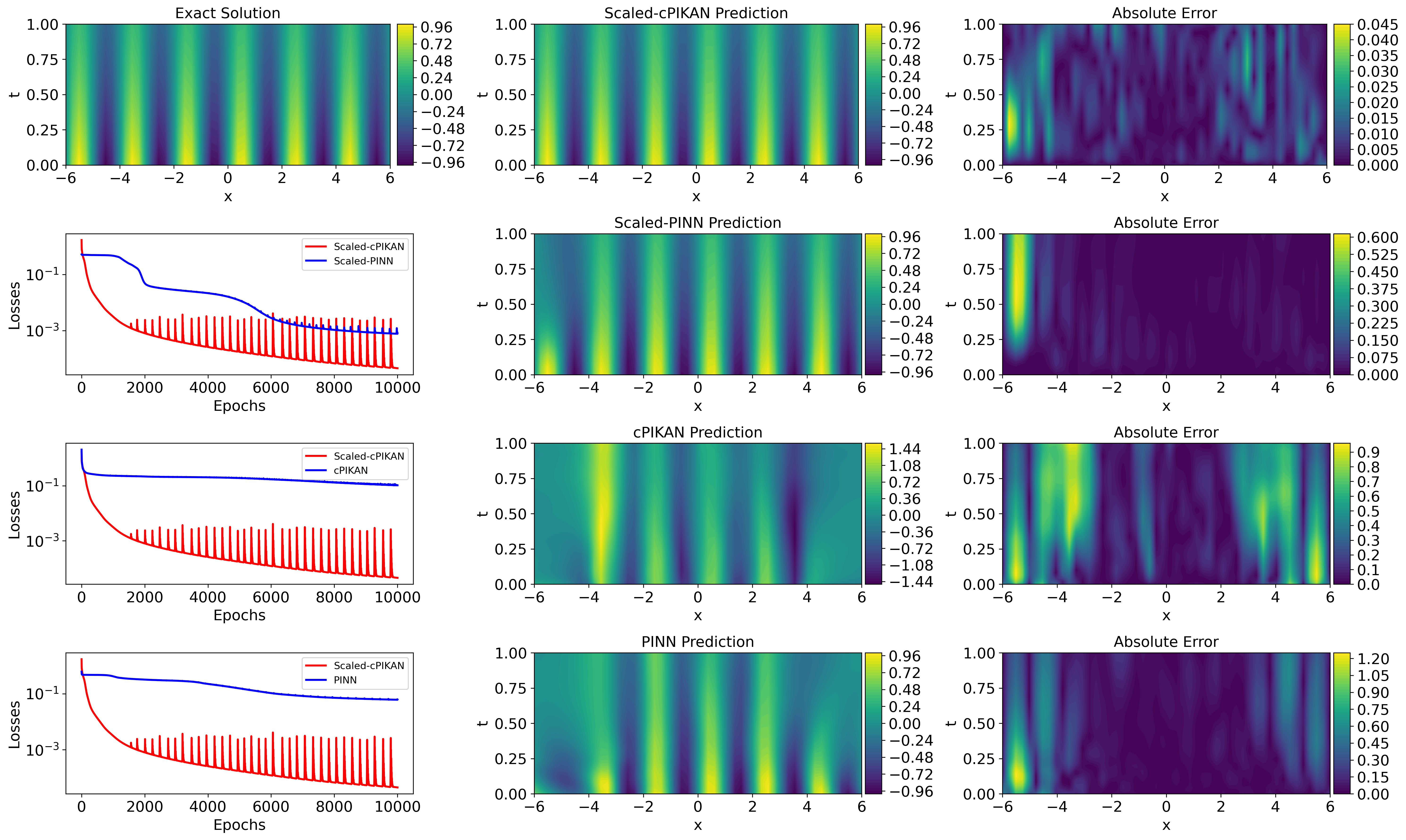}
    \caption{Comparison of the solutions predicted for the diffusion equation (Example \ref{Exam.DiffEqu}) in \([-6, 6] \times [0, 1]\). Each figure displays the outcomes of Scaled-cPIKAN, Scaled-PINN, cPIKAN, and PINN (from top to bottom). From left to right in each row: the loss function, the predicted solution \(u\), and the absolute error of the prediction are shown. The first plot in the top row shows the ground truth solution.}
    \label{fig:Heat6}
\end{figure}

We compare the results of these methods for various domain sizes in Figs.~\ref{fig:Heat2}–\ref{fig:Heat6}. Each figure presents the solution \(u\) over the domain \([-M, M] \times [0, 1]\), as predicted by Scaled-cPIKAN, Scaled-PINN, cPIKAN, and PINN, alongside the exact (ground truth) solution for comparison. Additionally, for enhanced comparison, the graphs for loss function versus epoch and absolute error are demonstrated. 
In Fig.~\ref{fig:Heat2}, it  can be seen that  the Scaled-cPIKAN model provides a prediction with a minimum absolute error, while capturing the fine details in the solution. The maximum absolute error in Scaled-cPIKAN is half that in vanilla cPIKAN, while the loss function values of the latter are comparable with a lower loss for Scaled-cPIKAN. Moreover, the Scaled-cPIKAN model achieves a maximum absolute error that is threefold smaller than Scaled-PINN and four times smaller compared to PINN, while consistently exhibiting a lower value for the loss function.\\
In Fig.~\ref{fig:Heat4}, it is evident that, with increasing domain size, the performance of scaled methods improves even significantly. Specifically, the maximum absolute error of Scaled-cPIKAN is half that of Scaled-PINN. Furthermore, the maximum absolute error of Scaled-cPIKAN is twelve times smaller than that of cPIKAN, highlighting a significant improvement in accuracy.  Compared to PINN, the maximum absolute error of Scaled-cPIKAN is seven times smaller.
Moreover, after 10,000 iterations, Scaled-cPIKAN achieves a much lower loss ($4.3\times 10^{-5}$) compared to the other three methods.
Scaled-PINN, cPIKAN, and PINN give larger values of the loss function, close to \(1.1\times 10^{-4}\), \(4.2\times 10^{-3}\) and \(1.0 \times 10^{-3}\), respectively.
In terms of the rate of convergence, it can also be seen that our proposed model, Scaled-cPIKAN, is faster than any other method discussed in this study. \\
In Fig.~\ref{fig:Heat6}, when the spatial domain is stretched to \([-6,6]\), both PINN and cPIKAN failed to approximate the solution, as their performance dramatically degrades with increasing domain size. In contrast, Scaled-cPIKAN provided a very accurate approximation, and its maximum absolute error is 26 times smaller than PINN, 20 times smaller than cPIKAN, and 13 times smaller than Scaled-PINN. 
The loss function value of scaled-cPIKAN is close to \(4.5 \times 10^{-5}\), which is significantly improved compared to scaled-PINN (\(7.9 \times 10^{-4}\)), cPIKAN (\(1.0 \times 10^{-1}\)) and PINN (\(6.1 \times 10^{-2}\)). These results highlight the robustness of Scaled-cPIKAN in handling extended spatial domains.

From Figs.~\ref{fig:Heat2}–\ref{fig:Heat6}, it can be concluded that the proposed scaled methods produce significantly more accurate solutions than those provided by the vanilla methods, especially for extended spatial domains involving the oscillatory behavior of the solution. 
Furthermore, when evaluating the values of the loss function for Scaled-cPIKAN against those derived from other methods, as depicted in the left column of Figs.~\ref{fig:Heat2}–\ref{fig:Heat6}, it is evident that Scaled-cPIKAN demonstrates lower loss values. Notably, the rate at which the loss in Scaled-cPIKAN decreases with epochs is higher, indicating a faster convergence and enhanced learning efficiency compared to other models, especially cPIKAN and PINN. These results show better performance of the scaled approach, especially when dealing with problems with high-frequency solutions over extended spatial domains.
% % 


\subsection{ Helmholtz Equation }\label{Exam.HelH}
As the second example, we aim to test the performance of the proposed model to solve the two-dimensional Helmholtz equation with homogeneous Dirichlet boundary conditions,
\begin{equation}
\begin{array}{ll}
u_{xx}(x,y) + u_{yy}(x,y) + \kappa^2 u(x,y) = f(x, y), & (x, y) \in \Omega,\\[3mm]
u(x,y) = 0, &(x, y) \in \partial \Omega ,
\end{array} 
\end{equation}
where the computational domain is defined as \(\Omega = [-M_x, M_x] \times [-M_y, M_y]\).
The forcing term is given by,
\begin{equation}
    f(x, y) = \left(\kappa^2 - \left(a_1^2 + a_2^2\right) \pi^2\right) \sin(a_1 \pi x) \sin(a_2 \pi y),
\end{equation}
where \(\kappa > 0\) is the wave number and \(a_1, a_2 \in \mathbb{R}\) are parameters controlling the frequency of the oscillations in the ground truth solution.
The ground truth solution for this example reads as, 
\begin{equation}
    u(x, y) = \sin(a_1 \pi x) \sin(a_2 \pi y),
\end{equation}
which satisfies the boundary conditions \(u(x, y) = 0\) on \(\partial \Omega\), since \(\sin(a_1 \pi x)\) and \(\sin(a_2 \pi y)\) vanish at the boundaries of \(\Omega\).
Changing the parameters \(\kappa, a_1,\) and \(a_2\) tune the problem to exhibit various levels of complexity, i.e., higher oscillatory behavior versus smoother solutions.
In this example, we solve the problem in a forward scheme without the presence of an initial condition. Consequently, the coefficients for the initial and measurement data terms in the data loss function, \(\lambda_{\text{init}}\) and \(\lambda_{\text{meas}}\), are set to zero in Eq.~\eqref{Eq.SLossData}.
We solve this example with two different settings on the domain \(\Omega = [-4, 4]^2\). 

Table \ref{Tab.HelmholtzSetting} presents the networks' configurations and the associated relative \(\mathcal{L}^2\) errors of the solutions predicted for the Helmholtz equation under two different PDE settings. The first case has \((a_1, a_2, \kappa) = (0.25, 1.0, 1.0)\) with equal weightings of the residual and boundary condition losses \((\lambda_{\text{res}}, \lambda_{\text{data}}) = (0.5, 0.5)\), and sampling distributions of 2000 domain points and 800 boundary points, $(N_{\text{res}}, N_{\text{data}}) = (2000, 800)$. In this case, Scaled-cPIKAN shows the best accuracy with the relative \(\mathcal{L}^2\) error of \(5.5 \times 10^{-2}\), outperforming Scaled-PINN \(1.6 \times 10^{-1}\), as well as  cPIKAN \(1.3 \times 10^0\) and PINN \(2.0 \times 10^{-1}\).\\
We further consider a solution with a higher frequency and choose $(a_1, a_2, \kappa) = (1.0, 1.0, 1.0)$, while also using a strongly weighted boundary loss $(\lambda_{\text{res}}, \lambda_{\text{data}}) = (0.01, 1.0)$ and sampling distributions of   $N_{\text{res}}=N_{\text{data}}=4000$. In this case, the smallest error is again attained by Scaled-cPIKAN. 
Comparing Scaled-cPIKAN with Scaled-PINN, Scaled-cPIKAN produces a much lower relative $\mathcal{L}^2$ error, raising the accuracy by roughly 98.4\% compared with Scaled-PINN. Also, compared with other methods, Scaled-cPIKAN decreases the relative $\mathcal{L}^2$ error by about 93.0\% to cPIKAN and 96.5\% to PINN. These results indicate the effectiveness of scaling transformations in enhancing the accuracy of neural network-based solvers, especially for challenging problems that have the oscillatory behavior of the solution.
% % %%%=============================================================
\begin{table}[!h]
\centering
\caption{\label{Tab.HelmholtzSetting} 
Network configurations, number of parameters, and associated relative \(\mathcal{L}^2\) errors  for solving the Helmholtz equation (Example \ref{Exam.HelH}) on \(\Omega = [-4, 4]^2\) under two distinct scenarios.}
\renewcommand{\arraystretch}{1.2}
\setlength{\tabcolsep}{12pt}
\begin{tabular}{c|cc|c}
\toprule
\multicolumn{4}{l}{$(a_1, a_2, \kappa) = (0.25, 1.0, 1.0)$,}\\
\multicolumn{4}{l}{$(\lambda_{\text{res}}, \lambda_{\text{data}}) = (0.5, 0.5)$, and $(N_{\text{res}}, N_{\text{data}}) = (2000, 800)$} \\
\midrule
% &  &  & \multicolumn{3}{c}{Relative $\mathcal{L}^2$ Error} \\
Method & $(N_l, N_n, k)$ & $\vert\boldsymbol{\theta} \vert$ & Relative $\mathcal{L}^2$ Error\\
\midrule
Scaled-cPIKAN & (4, 8, 3) & 864 &$\mathbf{5.5\times10^{-2}}$ \\
Scaled-PINN & (4, 16, -) & 960 & $1.6\times10^{-1}$\\
cPIKAN & (4, 8, 3) & 864 &$1.3\times 10^0$\\
PINN & (4, 16, -) & 960 & $2.0\times10^{-1}$\\
\midrule
\midrule
\multicolumn{4}{l}{$(a_1, a_2, \kappa) = (1.0, 1.0, 1.0)$,} \\
\multicolumn{4}{l}{$(\lambda_{\text{res}}, \lambda_{\text{data}}) = (0.01, 1.0)$, and $(N_{\text{res}}, N_{\text{data}}) = (4000, 4000)$} \\
\midrule
% &  &  & \multicolumn{3}{c}{Relative $\mathcal{L}^2$ Error} \\
Method & $(L, N, D)$ & $\vert\boldsymbol{\theta} \vert$ & Relative $\mathcal{L}^2$ Error\\
\midrule
Scaled-cPIKAN & (4, 25, 3) & 7800 &$\mathbf{3.2\times10^{-2}}$ \\
Scaled-PINN & (4, 50, -) & 7850 & $2.0\times 10^0$\\
cPIKAN & (4, 25, 3) & 7800&$4.6\times10^{-1}$\\
PINN & (4, 50, -) & 7850 & $9.1\times10^{-1}$\\

\bottomrule
\end{tabular}
\end{table}
% % % ===========================================

\begin{figure}[!h]
    \centering
\includegraphics[width=0.95\linewidth]{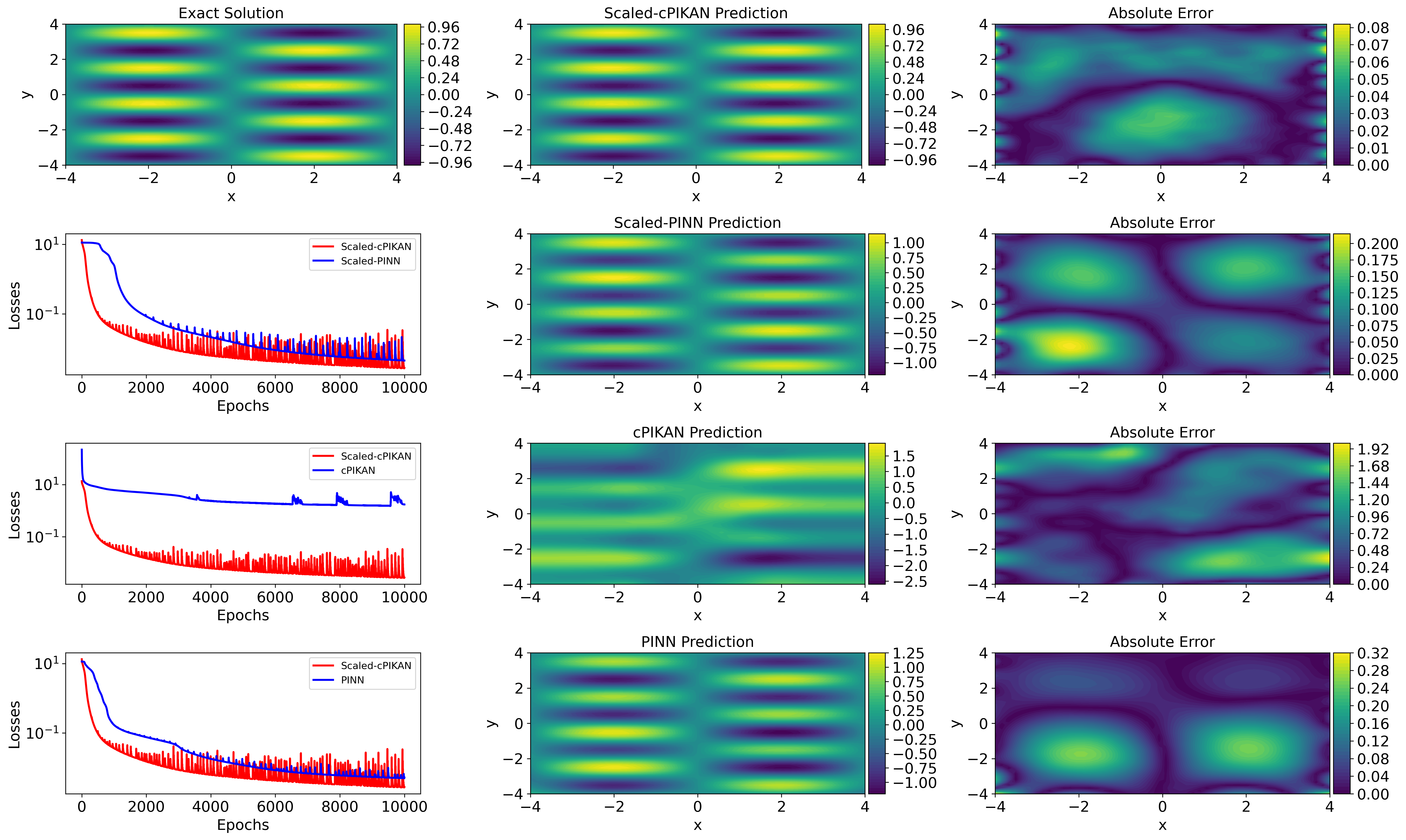}
    \caption{Prediction results for the Helmholtz equation (Example \ref{Exam.HelH}) in the first scenario with parameters \((a_1, a_2, \kappa) = (0.25, 1.0, 1.0)\) on \(\Omega = [-4,4]^2\). For each method (Scaled-cPIKAN, Scaled-PINN, cPIKAN, and PINN, shown sequentially from top to bottom), the plots from left to right display the loss function, the predicted solution \(u\), and the absolute error \(|u_{\text{exact}} - u_{\text{predicted}}|\). The first plot in the top row shows the ground truth solution.}
    \label{fig:Helm_1}
\end{figure}

\begin{figure}[!h]
    \centering
\includegraphics[width=0.95\linewidth]{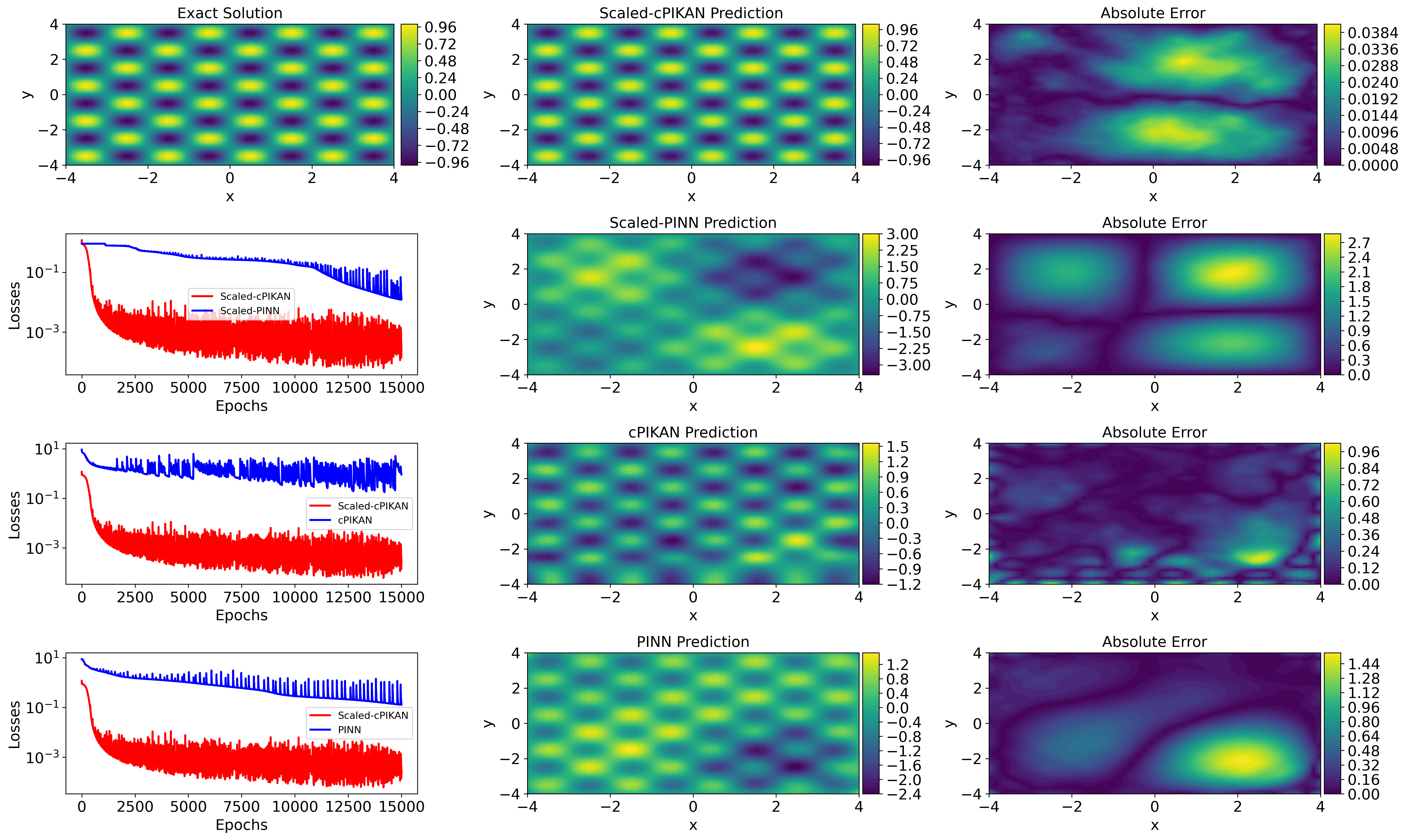}
    \caption{Prediction results for the Helmholtz equation (Example \ref{Exam.HelH}) in the second scenario with parameters \((a_1, a_2, \kappa) = (1.0, 1.0, 1.0)\) on \(\Omega = [-4,4]^2\). For each method (Scaled-cPIKAN, Scaled-PINN, cPIKAN, and PINN, shown sequentially from top to bottom), the plots from left to right display the loss function, the predicted solution \(u\), and the absolute error \(|u_{\text{exact}} - u_{\text{predicted}}|\). The first plot in the top row shows the ground truth solution.}
    \label{fig:Helm_2}
\end{figure}

The results of the four methods for the two simulation scenarios are shown in Figs.~\ref{fig:Helm_1}-\ref{fig:Helm_2}.  Figure \ref{fig:Helm_1} corresponds to the first scenario with parameters \((a_1, a_2, \kappa) = (0.25, 1.0, 1.0)\) in \(\Omega = [-4, 4]^2\). The comparison between Scaled-cPIKAN and Scaled-PINN indicates that Scaled-cPIKAN provides roughly twice as accurate results as Scaled-PINN by absolute error. Scaled-PINN and PINN provide similar approximations of \(u\), but Scaled-cPIKAN results in an error that is 24 times smaller compared to cPIKAN. 
After 10,000 iterations, Scaled-cPIKAN achieves a much lower loss ($2.7\times 10^{-3}$) compared to cPIKAN ($1.7 \times 10^{0}$). Similarly, Scaled-PINN slightly outperforms PINN ($4.5\times 10^{-3}$ vs. $5.1\times 10^{-3}$). These results highlight the effectiveness of scaling in both frameworks, especially when applied to KANs, i.e.,  Scaled-cPIKAN.

%FFFFFFFFFFFFFFFFFFFFFFFFFFFFFFFFFF
% % % %FFFFFFFFFFFFFFFFFFFFFFFFFFFFFFFFFF

Figure \ref{fig:Helm_2} presents the results for the second scenario with parameters $(a_1, a_2, \kappa) = (1.0, 1.0, 1.0)$ in the Helmholtz equation over the same domain $\Omega = [-4, 4]^2$. In this case, Scaled-cPIKAN consistently outperforms the other methods tested, including Scaled-PINN, cPIKAN, and PINN. Scaled-cPIKAN shows its strengths in capturing the oscillatory feature of the solution $u$, which is a challenging aspect of high-frequency solutions. 
Scaled-cPIKAN has a maximum absolute error that is 71, 25, and 37 times smaller than Scaled-PINN, cPIKAN, and PINN, respectively. 
The loss values further highlight the advantages of Scaled-cPIKAN, which achieves a minimum loss close to \(2.3 \times 10^{-4}\). In contrast, Scaled-PINN, cPIKAN, and PINN all reach a minimum loss of \(1.2 \times 10^{-2}\), \(9.1 \times 10^{-1}\), and \(1.3 \times 10^{-1}\), respectively, highlighting the superior performance of Scaled-cPIKAN in both accuracy and efficiency.
These results highlight the enhanced predictive capability and convergence efficiency of Scaled-cPIKAN for problems with different complexities, hence positioning it as a reliable and accurate solver for Helmholtz equations.
% % % 

% % % %FFFFFFFFFFFFFFFFFFFFFFFFFFFFFFFFFF

% \newpage
\subsection{Allen-Cahn Equation }\label{Exam.AC}
In the third example, we plan to test the efficacy of Scaled-cPIKAN for solving the Allen-Cahn equation widely used to model phase separation processes in multi-component systems. It is formulated as,
\begin{equation}\label{Eq.AC}
\begin{array}{cc}
     u_{t} - D \, u_{xx} + 5 \left(u^3 - u\right) = 0,& x \in [-M, M], \, t \in (0, T], \\[3mm]
     u(-M, t) = u(M, t) = -1, &\quad t \in (0, T], \\[3mm]
     u(x, 0) = (x/M)^2 \cos\left(\pi \:x/M\right), &\quad x \in [-M, M],
\end{array}
\end{equation}  
where \(u\) represents the state variable of interest, \(D\) is the diffusion coefficient, and the nonlinear term \(5(u^3 - u)\) models the local phase dynamics.
Solving this problem is challenging due to several reasons, (i) the nonlinear term \(5(u^3 - u)\) that leads to stiff behavior and potential difficulties in numerical stability; (ii) the diffusion term \(-D \, u_{xx}\) that requires fine spatial and temporal resolution to clearly capture the smooth transitions of the solution; and lastly (iii) the oscillatory behavior introduced by the initial condition that can propagate and interact with the nonlinear effects, yielding a highly sensitive evolution to numerical approximations. 
All of these factors combined call for robust numerical methods to ensure that solutions are both accurate and stable over the entire spatial domain and time.
Due to the mentioned challenges in solving the Allen-Cahn equation numerically, we address it in both supervised and unsupervised settings and report the results for each case separately.

In a supervised setup, the numbers of residual, boundary condition, initial condition, and measurement points are defined as \((N_{\text{res}}, N_{\text{bc}}, N_{\text{init}}, N_{\text{meas}}) = (2000, 400, 200, 200)\). The weighting factors for the residual and data loss terms are set to \((\lambda_{\text{res}}, \lambda_{\text{data}}) = (0.01, 1.0)\).  
For a fair comparison of the performance of the four methods, each network is again specified by its number of layers, neurons per layer, and polynomial degree, denoted as \((N_l, N_n, k)\), ensuring that they have approximately the same number of trainable parameters, \(\vert \boldsymbol{\theta} \vert\).  
Table \ref{Tab.ACSetting.Sup} summarizes the network configurations assigned to these methods to solve the Allen-Cahn equation with \(D = 10^{-4}\) over the domain \([-M, M]\times[0,1]\). The table also reports the relative \(\mathcal{L}^2\) errors for three different domain sizes: $[-2,2]\times[0,1]$, $[-4,4]\times[0,1]$,and $[-6,6]\times[0,1]$.  
The results show that Scaled-cPIKAN yields a higher accuracy for all sizes of domains. In contrast, Scaled-PINN has larger errors and does not provide an approximate solution for extended domains, e.g., (\(M=6\)). Both cPIKAN and PINN, the non-scaled versions, show severe degradation in performance as \(M\) increases, with both methods producing errors greater than \(0.5\) when \(M=6\). These results again highlight the advantages of scaled versions, especially Scaled-cPIKAN, in yielding accurate and stable solutions  on extended spatial domains.

% % %%%=============================================================
\begin{table}[!h]
\centering
\caption{\label{Tab.ACSetting.Sup} 
Network configurations, number of parameters, and associated relative \(\mathcal{L}^2\) errors for solving  the Allen-Cahn equation (Example \ref{Exam.AC}) in a supervised version with \(D = 10^{-4}\). The supervised training settings are \((N_{\text{res}}, N_{\text{bc}}, N_{\text{init}}, N_{\text{meas}}) = (2000, 400, 200, 200)\), with \((\lambda_{\text{res}}, \lambda_{\text{data}}) = (0.01, 1.0)\).}
\renewcommand{\arraystretch}{1.2}
\setlength{\tabcolsep}{12pt}
\begin{tabular}{c|cc|ccc}
\toprule
 &  &  & \multicolumn{3}{c}{Relative $\mathcal{L}^2$ Error} \\
Method & $(N_l, N_n, k)$ & $\vert \boldsymbol{\theta}\vert$ & $M=2$ & $M=4$ & $M=6$\\
\midrule
Scaled-cPIKAN & (4, 15, 3) & 2880 &$\mathbf{4.4\times 10^{-3}}$ & $\mathbf{7.8\times 10^{-3}}$& $\mathbf{6.9\times 10^{-3}}$\\
Scaled-PINN & (4, 30, -) & 2911 & $3.6\times 10^{-2}$& $2.4\times 10^{-2}$& $2.7\times 10^{-2}$\\
cPIKAN & (4, 15, 3) & 2880 &$6.1\times 10^{-2}$ & $5.9\times 10^{-1}$& $5.6\times 10^{-1}$\\
PINN &  (4, 30, -) & 2911 & $2.0\times 10^{-2}$&$5.2\times 10^{-1}$&$5.1\times 10^{-1}$\\
%\midrule
\bottomrule
\end{tabular}
\end{table}
% % %%%=============================================================

\begin{figure}
    \centering
    \includegraphics[width=0.95\linewidth]{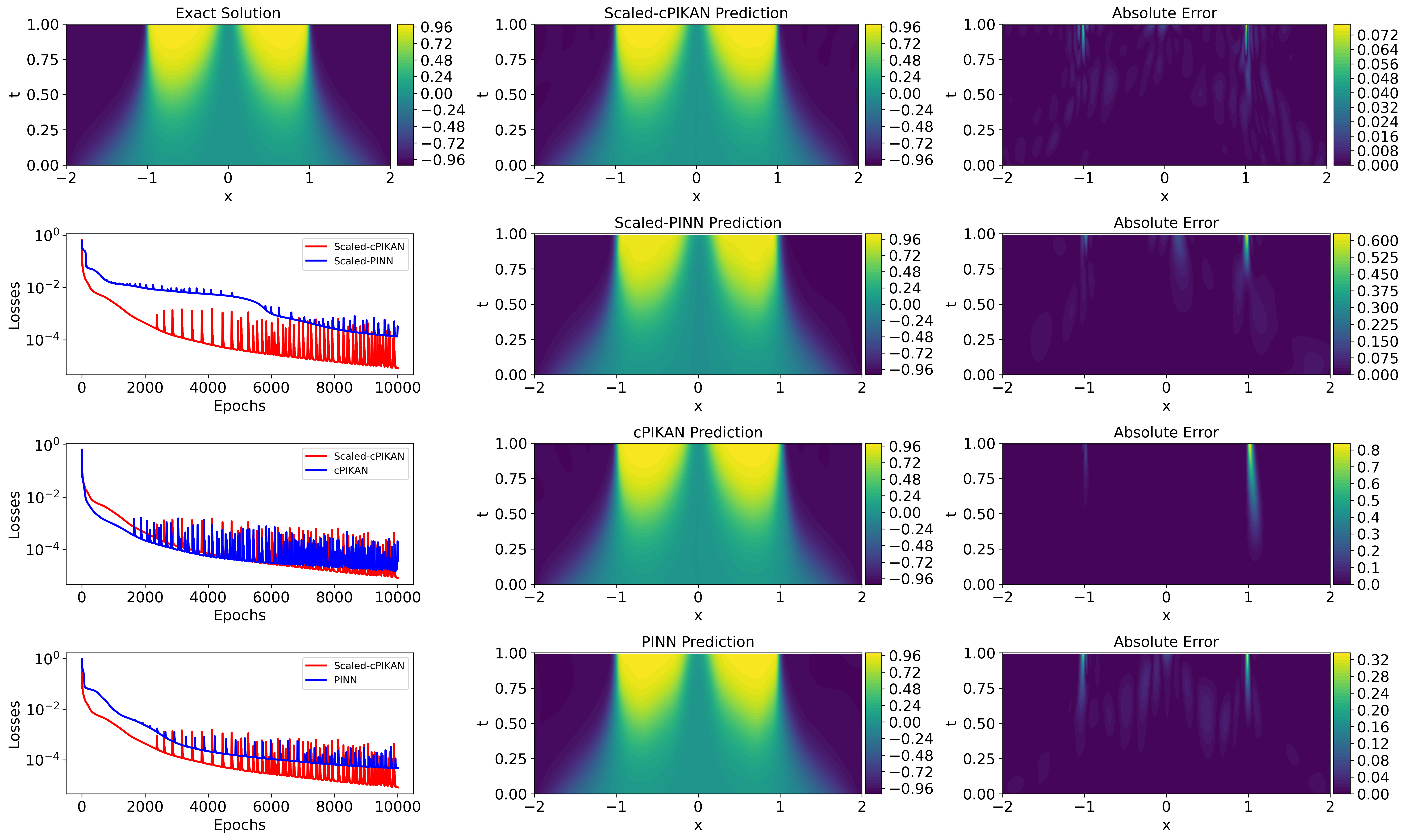}
    \caption{Comparison of the solutions predicted for the Allen-Cahn equation (Example \ref{Exam.AC}) in a supervised version on \([-2, 2]\times [0,1]\). Each figure displays the outcomes of Scaled-cPIKAN, Scaled-PINN, cPIKAN, and PINN (from top to bottom). From left to right in each row: the loss function, the predicted solution \(u\), and the absolute error of the prediction are shown. In the first row, the first plot shows the ground truth solution instead of the loss function for comparison.}
    \label{fig:AC2}
\end{figure}

\begin{figure}
    \centering
    \includegraphics[width=0.95\linewidth]{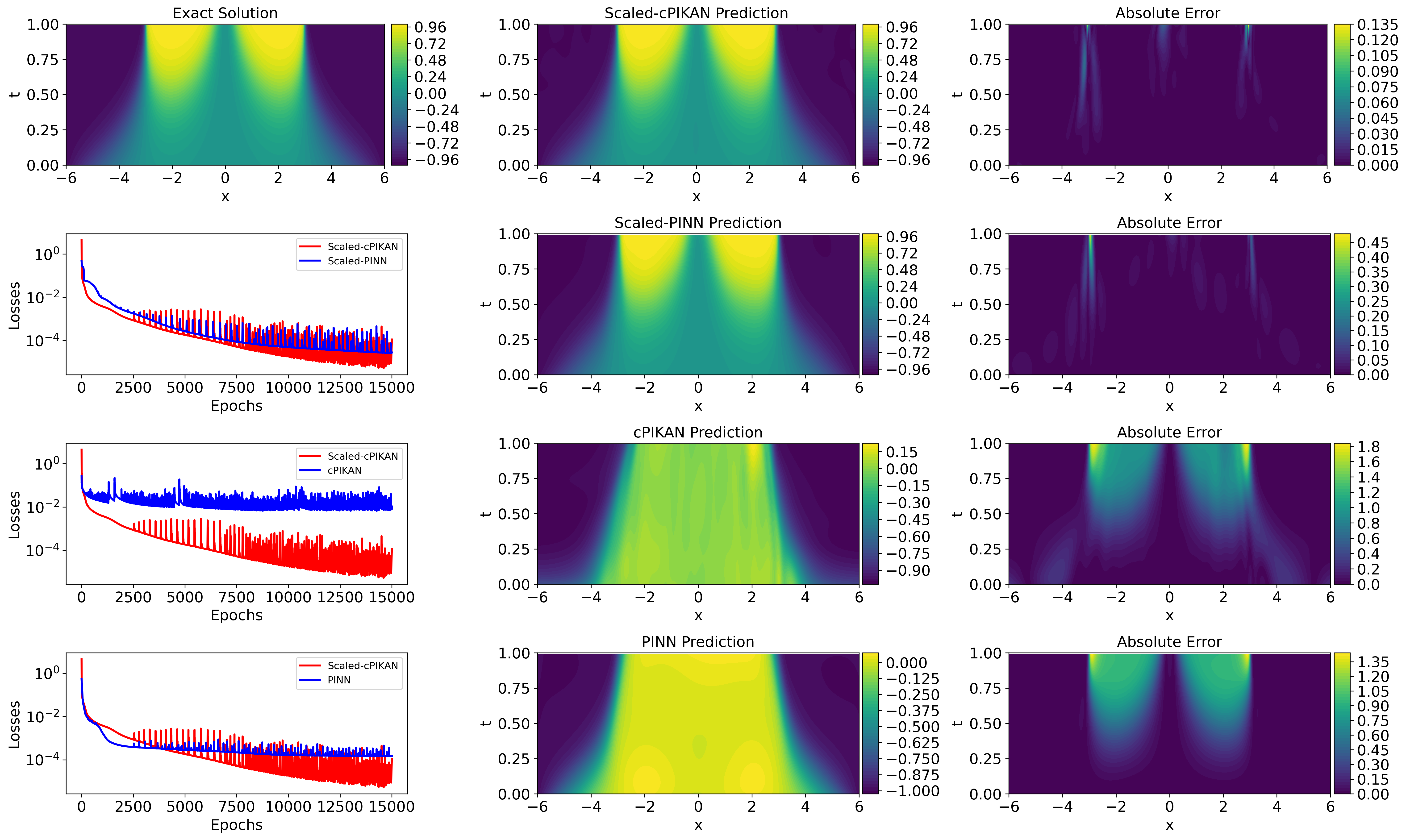}
    \caption{Comparison of the solutions predicted for the Allen-Cahn equation (Example \ref{Exam.AC}) in a supervised version on \([-6, 6]\times [0,1]\). Each figure displays the outcomes of Scaled-cPIKAN, Scaled-PINN, cPIKAN, and PINN (from top to bottom). From left to right in each row: the loss function, the predicted solution \(u\), and the absolute error of the prediction are shown. In the first row, the first plot shows the ground truth solution instead of the loss function for comparison.}
    \label{fig:AC6}
\end{figure}

Figures \ref{fig:AC2} and \ref{fig:AC6} represent the results of the Allen-Cahn equation in a supervised setting over different domains: $[-2,2]\times[0,1]$ and $[-6,6]\times[0,1]$. Figure~\ref{fig:AC2} ($M=2$) shows that while all four methods successfully predict the solution, Scaled-cPIKAN demonstrates superior performance. The maximum absolute error in Scaled-cPIKAN is notably eight times smaller than that in Scaled-PINN, 11 times smaller than that in cPIKAN, and four times smaller than that in PINN.  This reduction demonstrates the effect of the scaling transformation used in Scaled-cPIKAN.
Also, the convergence behavior of the loss function highlights the robustness of Scaled-cPIKAN. After 10,000 iterations, Scaled-cPIKAN achieves a loss as low as \( 8.2 \times 10^{-6} \), outperforming cPIKAN and PINN, which converge to losses of \( 4.4 \times 10^{-5} \) and \( 4.7 \times 10^{-5} \), respectively. In contrast, Scaled-PINN demonstrates reduced convergence speed, resulting in a final loss of \( 2.9 \times 10^{-4} \). This value is nearly two orders of magnitude higher compared to Scaled-cPIKAN.\\
In Fig.~\ref{fig:AC6} ($M=6$), as the domain size increases, cPIKAN and PINN show poor performance and cannot approximate the solution accurately. 
The maximum absolute error in Scaled-cPIKAN is significantly lower than the other methods, being three times lower than Scaled-PINN, eight times lower than cPIKAN, and ten times lower than PINN.  
The loss values also give a better view of the advantages of Scaled-cPIKAN. After 15,000 iterations, Scaled-cPIKAN achieves a loss that is roughly of the order of $10^{-6}$, e.g., drops to around $5.4 \times 10^{-6}$.
 In comparison, Scaled-PINN, cPIKAN, and PINN show a minimum loss of $2.8 \times 10^{-5}$, $1.0 \times 10^{-2}$, and $1.4 \times 10^{-4}$, respectively.
According to these results, both Scaled-cPIKAN and Scaled-PINN maintain better performance, with Scaled-cPIKAN consistently outperforming Scaled-PINN. 
Refer to Fig.~\ref{fig:AC4} in  \ref{App.Fig} for an illustration of the results when $M=4$, where Scaled-cPIKAN continues to exhibit greater accuracy than the other approaches.
% % % ===========================================
% % %FFFFFFFFFFFFFFFFFFFFFFFFFFFFFFFFFF

In the unsupervised context, by setting \(\lambda_{\text{meas}} = 0\) in Eq.~\eqref{Eq.SLossData}, we assess the effectiveness of the four methods in approximating the solution to the Allen-Cahn equation with \(D = 10^{-4}\) over the spatial domain \([-M, M]\) and temporal domain \([0,1]\). The simulations are carried out with \((\lambda_{\text{res}}, \lambda_{\text{data}}) = (0.01, 1.0)\), and the numbers of residual, boundary condition, and initial condition points are specified as \((N_{\text{res}}, N_{\text{bc}}, N_{\text{init}}) = (50000, 400, 200)\). The models are trained with a learning rate of \(10^{-4}\) over 150,000 epochs.
Table \ref{Tab.ACSetting.unSup} presents the relative \(\mathcal{L}^2\) errors for the different domain sizes \(M=2\), \(M=4\), and \(M=6\). According to the results reported in Table \ref{Tab.ACSetting.unSup}, we can see that Scaled-cPIKAN has the lowest error for all domain sizes. Specifically, for $M=2$, the error in Scaled-cPIKAN is over 90\% lower than that of Scaled-PINN and PINN, and approximately 40\% lower than cPIKAN.
As the value of $M$ increases, the gap between Scaled-cPIKAN's accuracy and  other methods is even more pronounced. For $M=4$, Scaled-cPIKAN maintains a significant advantage with an error reduction of almost 85\% compared to Scaled-PINN and PINN, and approximately 60\% better than cPIKAN.
By $M=6$, Scaled-cPIKAN continues to demonstrate superior performance, with the error being 90\% lower than Scaled-PINN and PINN, and around 90\% lower than cPIKAN. 
These results confirm the exceptional performance and robustness of Scaled-cPIKAN in producing very accurate approximations even with extended domain size, whereas all other approaches (in particular, PINN and Scaled-PINN) fail to provide reasonable approximations in the unsupervised setting.

% % %%%=============================================================
\begin{table}[!h]
\centering
\caption{\label{Tab.ACSetting.unSup} 
Network configurations, number of parameters, and associated relative \(\mathcal{L}^2\) errors for solving  the Allen-Cahn equation (Example \ref{Exam.AC}) in the unsupervised setting with \(D = 10^{-4}\). The training configuration includes \((N_{\text{res}}, N_{\text{bc}}, N_{\text{init}}) = (50000, 400, 200)\), \((\lambda_{\text{res}}, \lambda_{\text{data}}) = (0.01, 1.0)\), a learning rate of \(10^{-4}\), and 150,000 epochs.}
\renewcommand{\arraystretch}{1.2}
\setlength{\tabcolsep}{12pt}
\begin{tabular}{c|cc|ccc}
\toprule
&  &  & \multicolumn{3}{c}{Relative $\mathcal{L}^2$ Error} \\
Method & $(N_l, N_n, k)$ & $\vert \boldsymbol{\theta}\vert$  & $M=2$ & $M=4$ & $M=6$\\
\midrule
Scaled-cPIKAN & (4, 20, 5) & 7560 &$\mathbf{5.2\times 10^{-2}}$ & $\mathbf{7.8 \times10^{-2}}$& $\mathbf{5.9 \times10^{-2}}$\\
Scaled-PINN & (4, 50, -) & 7850 & $5.1\times 10^{-1}$& $5.1\times 10^{-1}$& $5.1\times 10^{-1}$\\
cPIKAN & (4, 20, 5) & 7560 &$8.9\times 10^{-2}$ & $2.4\times 10^{-1}$& $6.0 \times 10^{-1}$\\
PINN &  (4, 50, -) & 7850 & $5.1 \times 10^{-1}$&$5.1 \times 10^{-1}$&$5.1 \times 10^{-1}$\\
\bottomrule
\end{tabular}
\end{table}
% % % ===========================================

\begin{figure}
    \centering
    \includegraphics[width=0.95\linewidth]{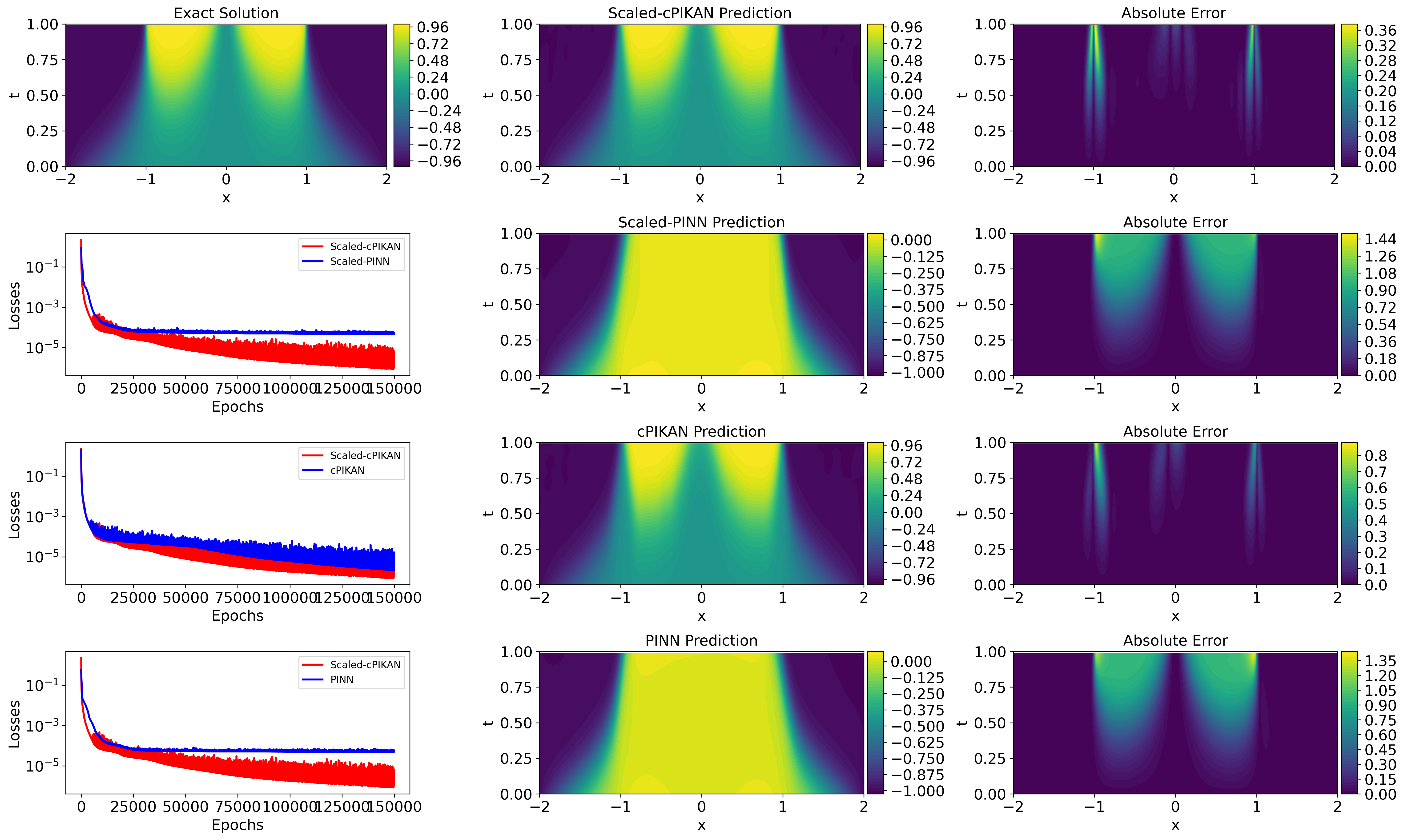}
    \caption{Comparison of the solutions predicted for the Allen-Cahn equation (Example \ref{Exam.AC}) in an unsupervised version on \([-2, 2]\times [0,1]\). Each figure displays the outcomes of Scaled-cPIKAN, Scaled-PINN, cPIKAN, and PINN (from top to bottom). From left to right in each row: the loss function, the predicted solution \(u\), and the absolute error of the prediction are shown. In the first row, the first plot shows the ground truth solution instead of the loss function for comparison.}
    \label{fig:unsAC2}
\end{figure}

\begin{figure}
    \centering
    \includegraphics[width=0.95\linewidth]{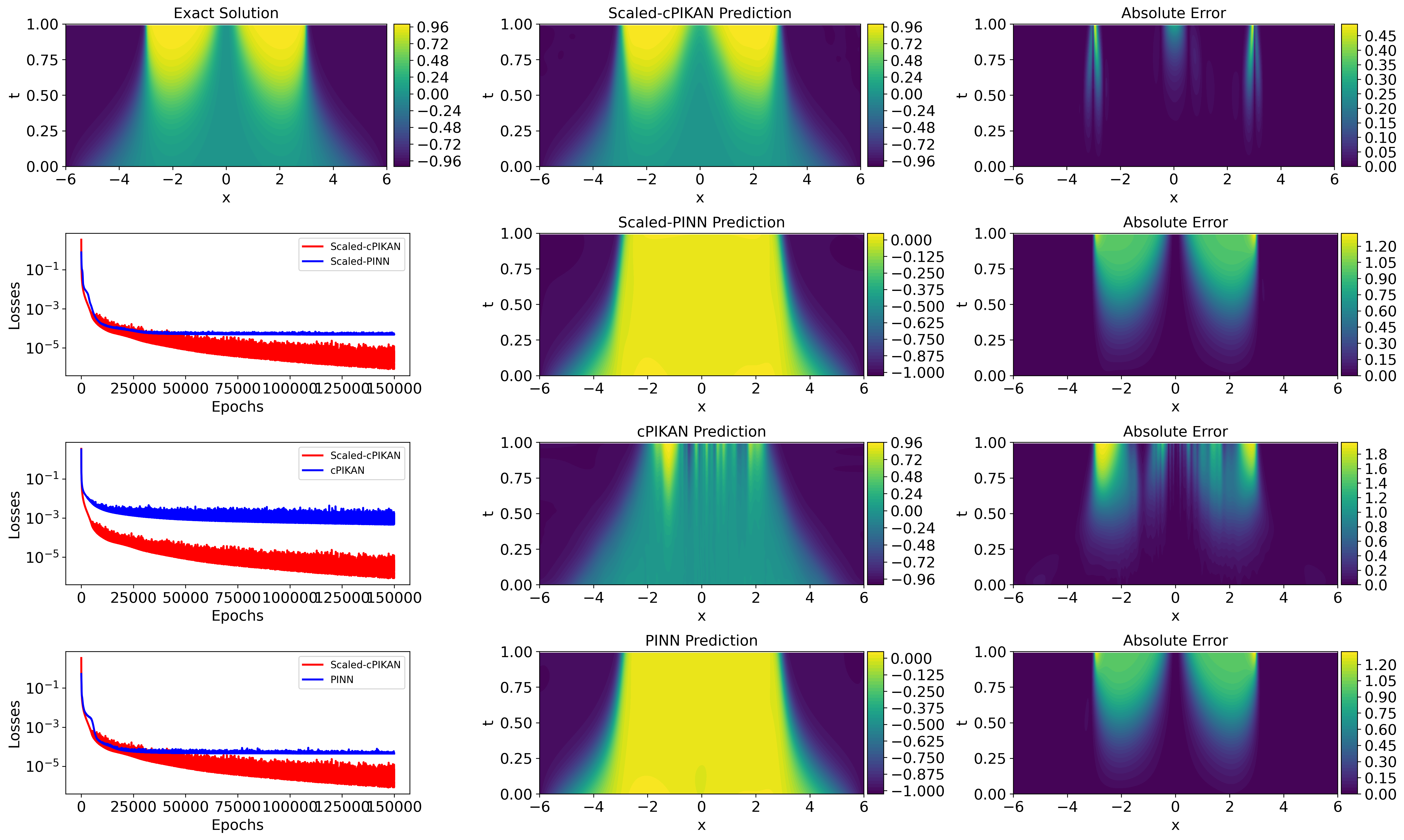}
    \caption{Comparison of the solutions predicted for the Allen-Cahn equation (Example \ref{Exam.AC}) in an unsupervised version on \([-6, 6]\times [0,1]\). Each figure displays the outcomes of Scaled-cPIKAN, Scaled-PINN, cPIKAN, and PINN (from top to bottom). From left to right in each row: the loss function, the predicted solution \(u\), and the absolute error of the prediction are shown. In the first row, the first plot shows the ground truth solution instead of the loss function for comparison.}
    \label{fig:unsAC6}
\end{figure}

The results of solving the Allen-Cahn equation in an unsupervised setting are presented in Figs.~\ref{fig:unsAC2} and \ref{fig:unsAC6} for spatial domains: \([-2, 2]\) and \([-6, 6]\), respectively. In Fig.~\ref{fig:unsAC2}, the cPIKAN-based methods provide a good approximation of the solution, with the Scaled-cPIKAN method achieving a maximum absolute error that is twice as small compared to the maximum error in cPIKAN. The loss function for both methods exhibits a similar behavior, steadily decreasing during the solution process and eventually converging to a value around \(10^{-6}\).
In contrast, the PINN-based methods fail to effectively approximate the solution.
Specifically, the maximum absolute error in Scaled-cPIKAN is 4 times smaller than both Scaled-PINN and PINN. 
A look at the nature of the loss function shows that Scaled-cPIKAN and cPIKAN achieve losses that are roughly on the order of $10^{-6}$, where Scaled-cPIKAN drops to $1.2 \times 10^{-6}$ and cPIKAN falls to $3.2 \times 10^{-6}$. On the other hand, Scaled-PINN and PINN losses fall within the order of $10^{-5}$, registering $4.8 \times 10^{-5}$ and $5.1 \times 10^{-5}$, respectively. 
\\
In Fig.~\ref{fig:unsAC6}, as the spatial domain expands, the Scaled-cPIKAN method provides a more satisfactory approximation of the solution over the entire time domain, while alternative methods fail to accurately approximate the function beyond time 0.5. The Scaled-cPIKAN method also demonstrates a maximum absolute error that is up to four times less than that of the cPIKAN method.
For Scaled-cPIKAN, the maximum absolute error is three times smaller than that of both the Scaled-PINN and PINN methods.
Moreover, the Scaled-cPIKAN method achieves a loss value on the order of \(10^{-6}\) in the end, while the cPIKAN method decreases to a value on the order of \(10^{-4}\). In contrast, both the Scaled-PINN and PINN methods decrease to loss values on the order of \(10^{-5}\).
For an illustration of $M=4$, see Fig.~\ref{fig:unsAC4} in \ref{App.Fig}, where Scaled-cPIKAN consistently exhibits superior accuracy over the other methods.

% % % ==========================================

% % %FFFFFFFFFFFFFFFFFFFFFFFFFFFFFFFFFF

% % %FFFFFFFFFFFFFFFFFFFFFFFFFFFFFFFFFF
\begin{figure}[!h]
\centering
\subfigure[]
{ \label{Pic.UNSAC06}\includegraphics[width=.8\textwidth]{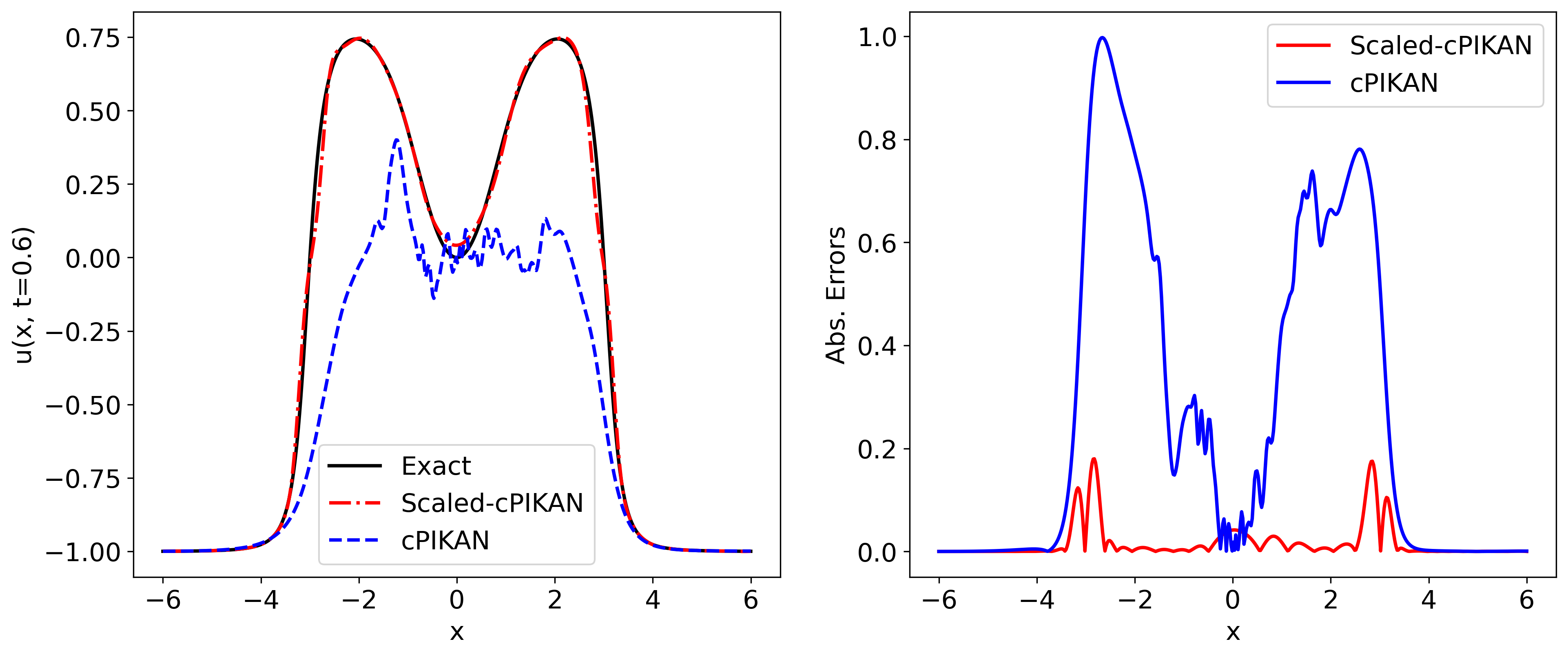}\hspace{.1cm}}
\subfigure[]
{ \label{Pic.UNSAC08}\includegraphics[width=.8\textwidth]{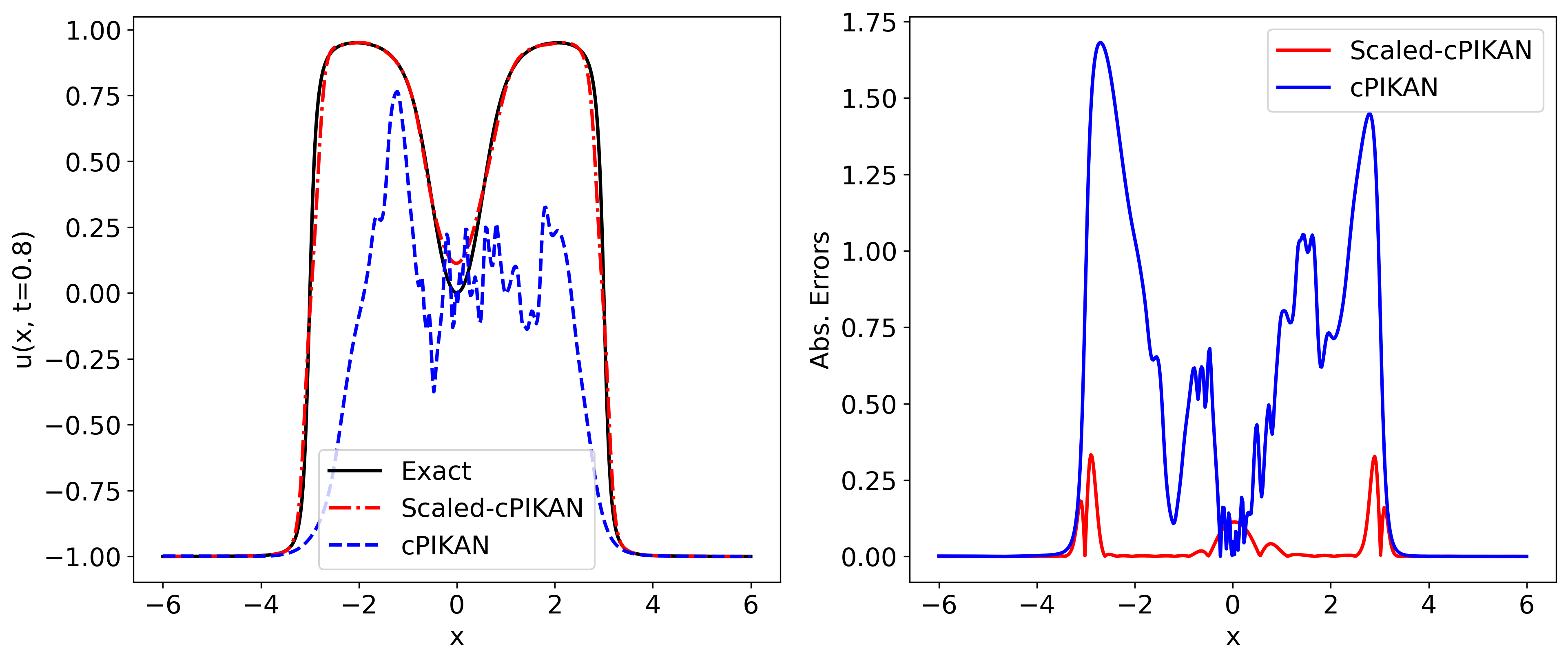}\hspace{.1cm}}
\subfigure[]
{ \label{Pic.UNSAC1}\includegraphics[width=.8\textwidth]{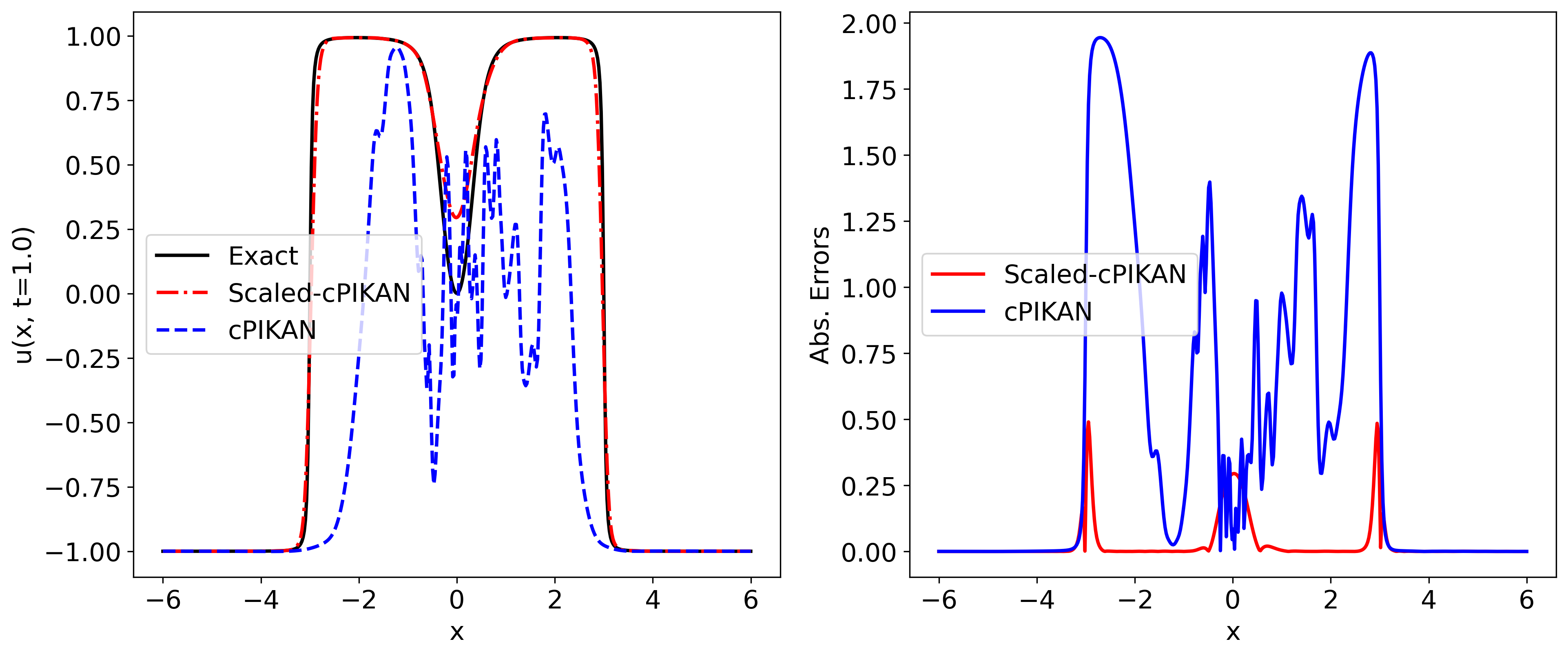}\hspace{.1cm}}
\caption{Comparing the performance of Scaled-cPIKAN and cPIKAN on solving the Allen-Cahn equation (Example \ref{Exam.AC}) in an unsupervised framework on the spatial domain \([-6,6]\) ] at three different time instances \(t = 0.6\), \(0.8\), and \(1.0\) (Panels (a)-(c), respectively). In each panel, the left plots show the predicted solutions \(u(x, t)\) together with the ground truth solution, while the right plots illustrate the corresponding absolute errors.} \label{Pic.UNSAC}%% label for entire figure
\end{figure}
% % 

The performance of Scaled-cPIKAN and cPIKAN in solving the Allen-Cahn equation under the unsupervised framework is compared in Figure \ref{Pic.UNSAC}. The results are presented for the spatial domain \([-6, 6]\) at three time instances \(t = 0.6\), \(0.8\), and \(1.0\) (Panels (a)–(c)).
In each panel, the left plot shows the predicted solutions \(u\) together with the ground truth solution, while the right plot depicts the corresponding absolute errors. It is clear that Scaled-cPIKAN outperforms cPIKAN by a large margin in the accuracy of its approximation and significantly reduces errors, particularly in capturing sharp transitions of the solution. These findings are in further concordance with the robustness and superiority of the scaling method in terms of accuracy and stability when applied over extended spatial domains and at various time points.

% % % ===========================================
\newpage
\subsection{Reaction-Diffusion Equation }\label{Exam.RD}
In the final example to test Scaled-cPIKAN's performance in both forward and inverse problems with noisy data, we consider the reaction-diffusion equation given by,  
\begin{equation}\label{Eq.RD}
    \begin{array}{ll}
D\, u_{xx} + \kappa\, \tanh(u) = f, & x \in \Omega, \\[3mm]
u(x) = \sin^3(6 x), & x \in \partial \Omega,
\end{array}
\end{equation} 
where $\Omega = [-M, M]$, $D$ represents the diffusion coefficient, and $\kappa$ is a reaction term coefficient, affecting the reaction dynamics. The ground truth solution to the problem is  \(u(x) = \sin^3(6 x)\). The source term \(f(x)\) is derived analytically by substituting the ground truth solution in Eq.~\eqref{Eq.RD}. This example presents a far-reaching test case for evaluating the accuracy and stability of our proposed methods. 
The second-order diffusion term \(D\: u_{xx}\) and the nonlinear reaction term \(\kappa\: \tanh(u)\) create complex interactions that make it challenging to approximate, as they involve sharp spatial changes and strong nonlinear effects.
We further demonstrate the strengths and limitations of our model in dealing with complex reaction-diffusion problems with non-trivial boundary conditions by comparing the performance of different approaches in this example. 

The parameters of the equation for the forward version are set to \(D = 10^{-2}\) and \(\kappa = 0.7\), and the domain \(\Omega = [-M, M]\) is tested for three sizes: \(M = 2\), \(4\), and \(6\). The models use $N_{\text{res}}=3000$ residual points with $\lambda_{\text{res}} = 1.0$, a learning rate of $10^{-4},$ and 20,000 epochs in training. For fairness in comparison, the neural networks are again set up to be similar in the number of parameters. Scaled-cPIKAN and cPIKAN had a configuration of \(N_l = 4\) hidden layers, \(N_n = 8\) neurons per layer, and \(k = 5\) Chebyshev polynomials, totaling 1248 parameters. Similarly, Scaled-PINN and PINN possess \(N_l = 4\) layers and \(N_n = 20\) neurons per layer, giving 1240 parameters.
For the forward version of the reaction-diffusion equation example, the network settings and the relative \(\mathcal{L}^2\) errors are summarized in Table \ref{Tab.RDSettingForward}. 
The results point to much better performance for the Scaled-cPIKAN and cPIKAN architectures in terms of accuracy across all domain sizes. 
For \(M = 2\), the relative \(\mathcal{L}^2\) errors of cPIKAN-based methods (Scaled-cPIKAN and cPIKAN) are comparable; however, that of Scaled-cPIKAN is 96.8\% less than that of Scaled-PINN, and 96.9\% less than that of PINN.
For \( M=6 \), only the Scaled-cPIKAN method retains its ability to accurately approximate the solution as the size of the domain increases.
The Relative \(\mathcal{L}^2\) Error of Scaled-cPIKAN is, on average, approximately 90\% lower than that of Scaled-PINN, cPIKAN, and PINN.

% % %%%=============================================================
\begin{table}[!h]
\centering
\caption{\label{Tab.RDSettingForward} 
Network configurations, number of parameters, and associated relative \(\mathcal{L}^2\) errors for solving  the reaction-diffusion equation (Example \ref{Exam.RD}) in the forward version on $\Omega=[-M, M]$. The parameters are \(D = 10^{-2}\), \(\kappa = 0.7\), and the training settings include \(N_{\text{res}} = 3000\) residual points, \(\lambda_{\text{res}} = 1.0\), a learning rate of \(10^{-4}\), and 20,000 training epochs.}
\renewcommand{\arraystretch}{1.2}
\setlength{\tabcolsep}{12pt}
\begin{tabular}{c|cc|ccc}
\toprule
% \multicolumn{6}{l}{Forward version:}  \\
% \multicolumn{6}{l}{} \\
% \midrule
 &  &  & \multicolumn{3}{c}{Relative $\mathcal{L}^2$ Error} \\
Method & $(N_l, N_n, k)$ & Num. Param. & $M=2$ & $M=4$ & $M=6$\\
\midrule
Scaled-cPIKAN & (4, 8, 5) & 1248 &$\mathbf{2.7\times 10^{-2}}$ & $\mathbf{2.0\times 10^{-1}}$& $\mathbf{9.6\times 10^{-2}}$\\
Scaled-PINN & (4, 20, -) & 1240 & $8.5\times 10^{-1}$& $9.2\times 10^{-1}$& $9.7\times 10^{-1}$\\
cPIKAN & (4, 8, 5) & 1248&$2.5\times 10^{-2}$ & $4.6\times 10^{-1}$& $1.0\times 10^0$\\
PINN &  (4, 20, -) & 1240 & $8.6\times 10^{-1}$&$1.1\times 10^0$&$1.0\times 10^0$\\
\bottomrule
\end{tabular}
\end{table}
% % % ===========================================

\begin{figure}
    \centering
    \includegraphics[width=0.95\linewidth]{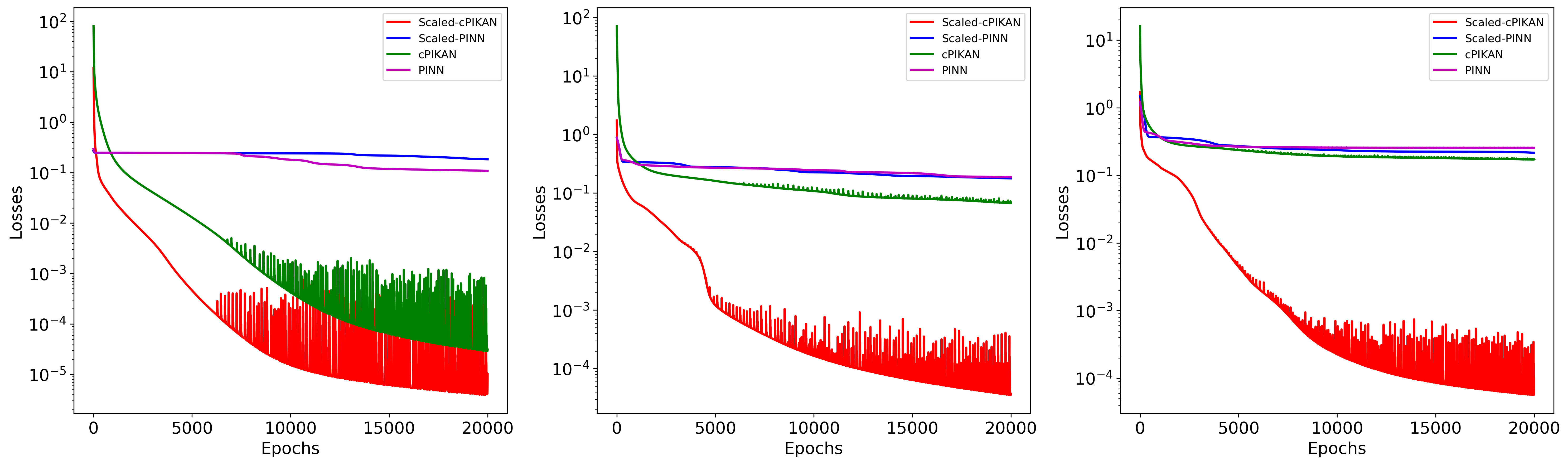}
    \caption{Training loss histories for the forward version of the reaction-diffusion equation (Example \ref{Exam.RD}) across three domain sizes: \(M = 2\), \(M = 4\), and \(M = 6\), displayed from left to right panels.}
    \label{fig:RD_Loss}
\end{figure}

Figure \ref{fig:RD_Loss} presents the training loss histories of the three different domain sizes for the forward problem. The proposed Scaled-cPIKAN has the minimum training loss with the fastest convergence.
In Fig.~\ref{fig:RD_Loss}, left panel (\( M=2 \)), the loss of the Scaled-cPIKAN method decreases very rapidly, reaching as low as \( 10^{-6} \). The cPIKAN method also shows a fast reduction, stabilizing around \( 10^{-5} \). In contrast, the two PINN-based methods oscillate around \( 10^{-1} \) after a few iterations without significant improvement.
In Fig.~\ref{fig:RD_Loss}, middle panel (\( M=4 \)), the Scaled-cPIKAN method continues to exhibit superior performance, with the loss rapidly decreasing to \( 10^{-5} \). However, the performance of the cPIKAN method weakens, with the loss only reducing to \( 10^{-2} \). The two PINN-based methods remain weak, oscillating around \( 10^{-1} \) after a few iterations, similar to the left panel.
In Fig.~\ref{fig:RD_Loss}, right panel (\( M=6 \)), only the Scaled-cPIKAN method performs well with the loss decreasing rapidly and consistently to \( 10^{-5} \). However, the other three methods stagnate, oscillating around \( 10^{-1} \) after a few iterations without a notable reduction in loss.\\
In addressing the forward formulation of the reaction-diffusion equation, it is observed that both PINN and Scaled-PINN approaches exhibit slower convergence rates and achieve higher final loss values compared to Scaled-cPIKAN. Consequently, these findings highlight the robustness of Scaled-cPIKAN and cPIKAN in managing the complex nonlinearities and sharp gradients inherent in the reaction-diffusion equation.

% % % ===========================================
%%%=============================================================

\begin{figure}
    \centering
    \includegraphics[width=0.95\linewidth]{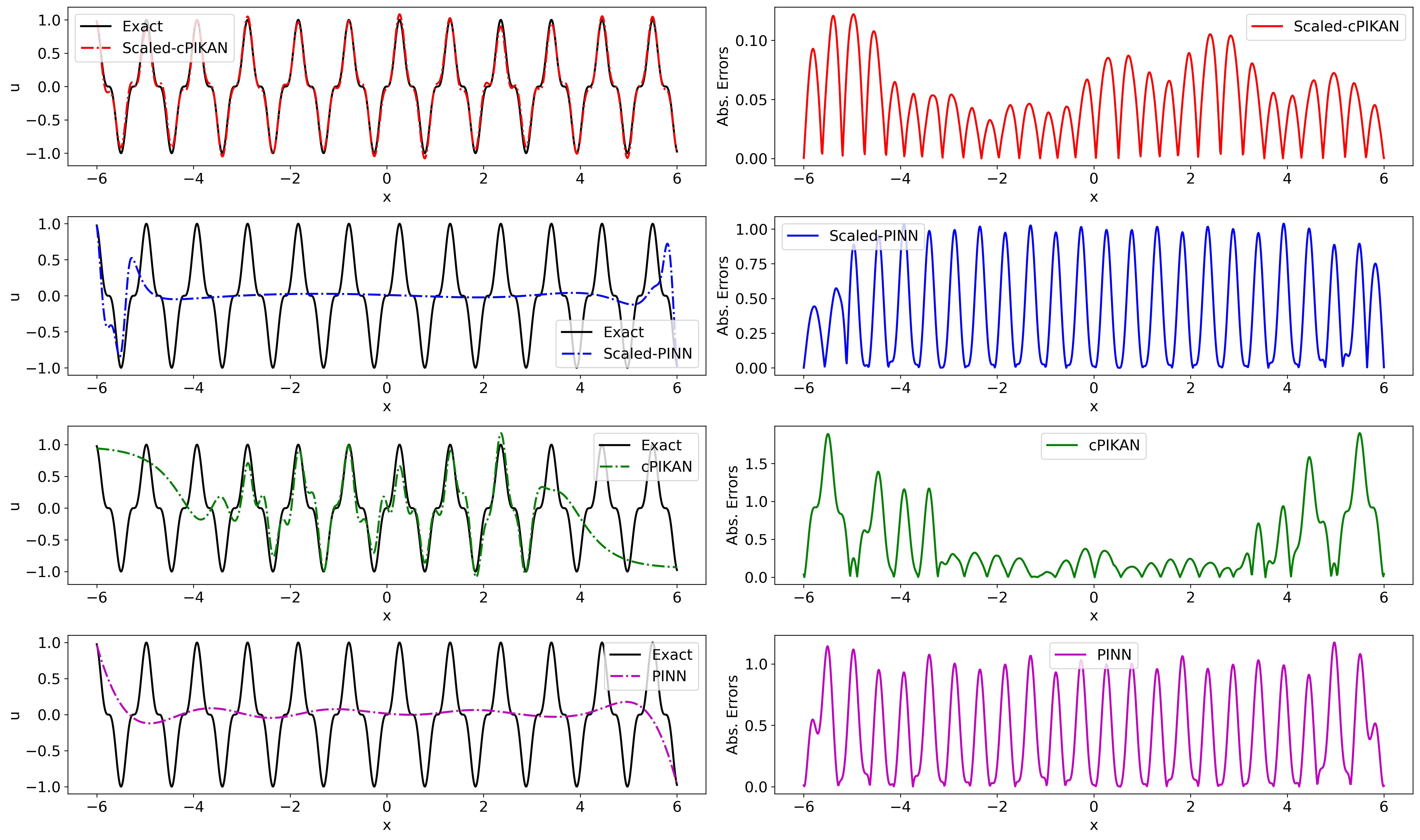}
    \caption{Forward solutions of the reaction-diffusion equation (Example \ref{Exam.RD}) for $M=6$ using Scaled-cPIKAN, Scaled-PINN, cPIKAN, and PINN methods (from top to bottom). Left Panels: comparison between the ground truth solution and predictions of the methods. Right Panels: absolute error distributions.}
    \label{fig:FRD6}
\end{figure}

In addition, Fig.~\ref{fig:FRD6} presents the forward solutions of the reaction-diffusion equation for domain sizes \(M = 6\), computed using four different methods, Scaled-cPIKAN, Scaled-PINN, cPIKAN, and PINN. 
In Fig.~\ref{fig:FRD6}, the Scaled-cPIKAN algorithm demonstrates the highest accuracy, achieving minimal errors in the extended domain \( [-6, 6] \) with a maximum absolute error of less than \( 0.12 \). In contrast, Scaled-PINN and PINN exhibit significant deviations from the ground truth solution, with maximum absolute errors below \( 1.00 \), making them less stable and accurate. Additionally, the cPIKAN method struggles to approximate the solution, particularly near the boundaries, resulting in a maximum absolute error of approximately \( 1.55 \). 
For an illustration of \( M=2 \) and \( M=4 \), refer to Figs.~\ref{fig:FRD2} and \ref{fig:FRD4} in  \ref{App.Fig}, where Scaled-cPIKAN and cPIKAN effectively capture the intricate features of the solution with minimal absolute error.
%
%%%%%%========================================================
% % % ===========================================

Next, we examine an inverse form of the reaction-diffusion equation given in Eq.~\eqref{Eq.RD} together with an estimation of the unknown parameter  \(\kappa\) and seek the neural networks approximations to solution \(u\). We also impose realistic noise on the measured data. Noisy observations of the solution and the source term are generated by adding random Gaussian noise with zero mean and prescribed standard deviations. 
The noisy data is defined as,  
\begin{equation}\label{Eq.NoiseU}
    \boldsymbol{u}^{\delta}_{\text{meas}} = \boldsymbol{u}_{\text{meas}} + \boldsymbol{\epsilon}_u, 
\end{equation}
and the noisy source term defined as,  
\begin{equation}\label{Eq.NoiseF}
    \boldsymbol{f}^{\delta} = \boldsymbol{f} + \boldsymbol{\epsilon}_f,
\end{equation}  
where \(\boldsymbol{\epsilon}_u\) and \(\boldsymbol{\epsilon}_f\) are vectors sampled from normal distributions with zero mean and variances \(\delta_u^2\) and \(\delta_f^2\), respectively.  
In Eq.~\eqref{Eq.NoiseU}, \(\boldsymbol{u}_{\text{meas}}\) is the vector representation of \(u_{\text{meas}}(\boldsymbol{\xi}_j^{\text{meas}}), j=1, \ldots, N_{\text{meas}}\) corresponding to the measurement values in the computational domain, and in Eq.~\eqref{Eq.NoiseF}, \(\boldsymbol{f}\) is the vector representation of \(f(\boldsymbol{\xi}_j^{\text{res}}), j=1, \ldots, N_{\text{res}}\) corresponding to the source values.  This approach investigates the robustness of our proposed model when dealing with noisy data in inverse problems, thus showcasing its capability for extracting reliable parameters even from noisy data.

The experiments are carried out in the domain \(\Omega = [-6, 6]\) using \((N_{\text{res}}, N_{\text{data}}) = (1200, 10)\) residual and data points, with equal weighting coefficients \((\lambda_{\text{res}}, \lambda_{\text{data}}) = (1.0, 1.0)\). The training is performed over 20,000 epochs with a learning rate of \(10^{-3}\). To ensure fairness, all methods utilize a similar network architecture with comparable parameter counts. In addition, three noise levels are considered: \((\delta_u, \delta_f) = (0.00, 0.00)\), \((\delta_u, \delta_f) = (0.05, 0.00)\), and \((\delta_u, \delta_f) = (0.05, 0.05)\). The corresponding network settings and relative \(\mathcal{L}^2\) errors for the inverse reaction-diffusion problem are summarized in Table \ref{Tab.RDSettingINV}.
The relative \(\mathcal{L}^2\) errors for solution \(u\), source term \(f\), and parameter \(\kappa\) are also reported in this table.

% % % ===========================================
\begin{table}[!h]
\centering
\caption{\label{Tab.RDSettingINV} 
Network configurations and relative \(\mathcal{L}^2\) errors for the reaction-diffusion equation (Example \ref{Exam.RD}) in the inverse version on $\Omega=[-6, 6]$. The training settings include \((N_{\text{res}}, N_{\text{data}}) = (1200, 10)\) residual points, \((\lambda_{\text{res}}, \lambda_{\text{data}}) = (1.0, 1.0)\), a learning rate of \(10^{-3}\), and 20,000 training epochs.}
\renewcommand{\arraystretch}{1.2}
\setlength{\tabcolsep}{12pt}
\begin{tabular}{c|cc|ccc}
\toprule
\multicolumn{6}{l}{$\delta_u = 0.00$ and $\delta_f = 0.00$}  \\
% \multicolumn{6}{l}{} \\
\midrule
 &  &  & \multicolumn{3}{c}{Relative $\mathcal{L}^2$ Error} \\
Method & $(N_l, N_n, k)$ & $\vert \boldsymbol{\theta}\vert$ & $u$ & $f$ & $\kappa$\\
\midrule
Scaled-cPIKAN & (4, 8, 5) & 1248 &$\mathbf{1.3\times 10^{-2}}$ & $\mathbf{9.9\times 10^{-3}}$& $\mathbf{6.8\times 10^{-4}}$\\
Scaled-PINN & (4, 20, -) & 1240 & $2.1\times 10^0$& $7.1\times 10^{-1}$& $5.7\times 10^{-1}$\\
cPIKAN & (4, 8, 5) & 1248&$1.9\times 10^0$ & $6.6\times 10^{-1}$& $7.6\times 10^{-1}$\\
PINN &  (4, 20, -) & 1240 & $3.1\times 10^0$&$6.8\times 10^{-1}$&$9.8\times 10^{-1}$\\
\midrule
\midrule
\multicolumn{6}{l}{$\delta_u = 0.05$ and $\delta_f = 0.00$}  \\
% \multicolumn{6}{l}{} \\
\midrule
 &  &  & \multicolumn{3}{c}{Relative $\mathcal{L}^2$ Error} \\
Method & $(N_l, N_n, k)$ & $\vert \boldsymbol{\theta}\vert$ & $u$ & $f$ & $\kappa$\\
\midrule
Scaled-cPIKAN & (4, 8, 5) & 1248 &$\mathbf{6.5\times 10^{-2}}$ & $\mathbf{1.6\times 10^{-2}}$& $\mathbf{3.8\times 10^{-4}}$\\
Scaled-PINN & (4, 20, -) & 1240 & $4.8\times 10^0$& $6.0\times 10^{-1}$& $9.2\times 10^{-1}$\\
cPIKAN & (4, 8, 5) & 1248&$2.0\times 10^0$ & $7.4\times 10^{-1}$& $8.3\times 10^{-1}$\\
PINN &  (4, 20, -) & 1240 & $3.2\times 10^0$&$5.9\times 10^{-1}$&$7.9\times 10^{-1}$\\
\midrule
\midrule
\multicolumn{6}{l}{$\delta_u = 0.05$ and $\delta_f = 0.05$}  \\
% \multicolumn{6}{l}{} \\
\midrule
 &  &  & \multicolumn{3}{c}{Relative $\mathcal{L}^2$ Error} \\
Method & $(N_l, N_n, k)$ & $\vert \boldsymbol{\theta}\vert$ & $u$ & $f$ & $\kappa$\\
\midrule
Scaled-cPIKAN & (4, 8, 5) & 1248 &$\mathbf{1.1\times 10^{-1}}$ & $\mathbf{5.3\times 10^{-2}}$& $\mathbf{4.6\times 10^{-3}}$\\
Scaled-PINN & (4, 20, -) & 1240 & $2.9\times 10^0$& $7.1\times 10^{-1}$& $7.8\times 10^{-1}$\\
cPIKAN & (4, 8, 5) & 1248&$1.9\times 10^0$ & $7.0\times 10^{-1}$& $8.7\times 10^{-1}$\\
PINN &  (4, 20, -) & 1240 & $2.6\times 10^0$&$4.9\times 10^{-1}$&$6.4\times 10^{-1}$\\
\bottomrule
\end{tabular}
\end{table}
% % % ===========================================

\begin{figure}[!h]
    \centering
    \includegraphics[width=0.95\linewidth]{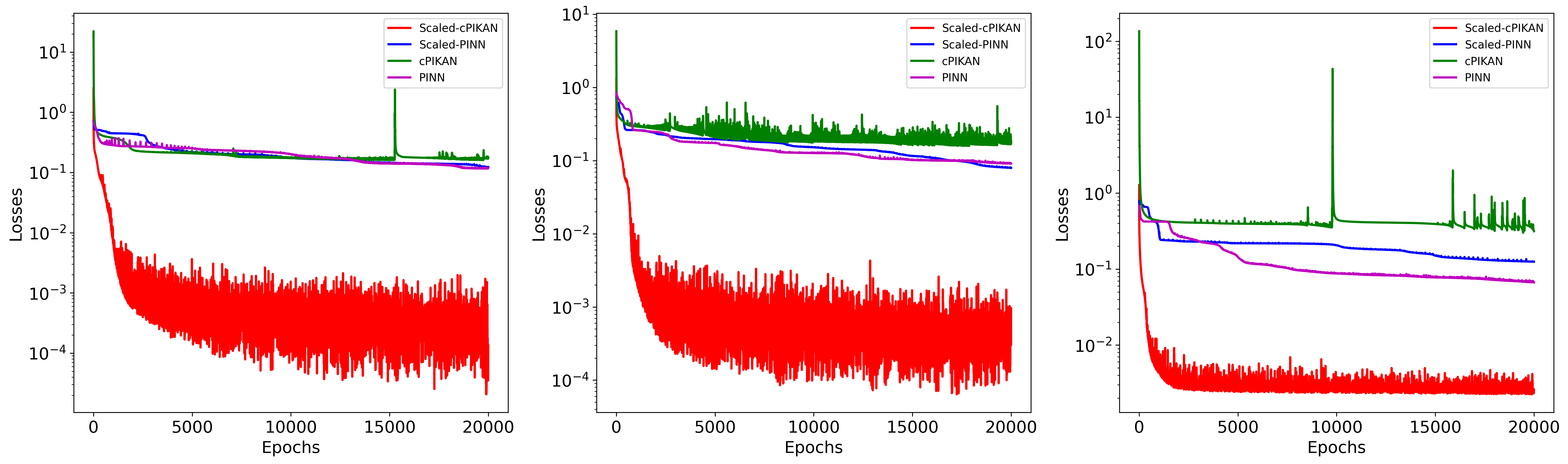}
    \caption{Training loss histories for the inverse version of the reaction-diffusion equation (Example \ref{Exam.RD}) across different noise levels. 
Left Panel: \(\delta_u = 0.00\) and \(\delta_f = 0.00\),
Middle Panel: \(\delta_u = 0.05\) and \(\delta_f = 0.00\),
Right Panel: \(\delta_u = 0.05\) and \(\delta_f = 0.05\).}
    \label{fig:INV_RD_Loss}
\end{figure}
% % % ===========================================
In Table \ref{Tab.RDSettingINV}, in the noiseless case (\(\delta_u = 0.00\) and \(\delta_f = 0.00\)), the Scaled-cPIKAN method achieves the most accurate results among all methods.
Specifically, it estimates \(\kappa\) with relative \(\mathcal{L}^2\) errors that are, on average, about 99\% lower than those of Scaled-PINN, cPIKAN, and PINN. Similarly, it approximates \(f\) with relative errors reduced by approximately 98\% and estimates \(u\) with relative errors that are around 99\% lower than the other methods.
When noise is introduced to the data (\(\delta_u = 0.05\) and \(\delta_f = 0.00\)), the Scaled-cPIKAN method maintains its superior performance, with \(\kappa\) estimates exhibiting relative \(\mathcal{L}^2\) errors that are, on average, about 99\% lower than those of Scaled-PINN, cPIKAN, and PINN. 
Its estimates of \(f\) show relative error reductions of 97.3\%, 97.8\%, and 97.3\%, and \(u\) estimates are 98.6\%, 96.8\%, and 97.9\% lower in comparison to other methods.
Even under the most challenging conditions (\(\delta_u = 0.05\) and \(\delta_f = 0.05\)), Scaled-cPIKAN continues to outperform all other methods. 
Our proposed method estimates \(\kappa\), \(f\), and \(u\) with average relative error reductions of approximately 99\%, 91\%, and 95\%, respectively, compared to Scaled-PINN, cPIKAN, and PINN.
These results, summarized in the table, highlight the robustness and accuracy of the Scaled-cPIKAN method, especially in noisy and noiseless scenarios.

% % % ===========================================
\begin{figure}[!h]
    \centering
    \includegraphics[width=0.95\linewidth]{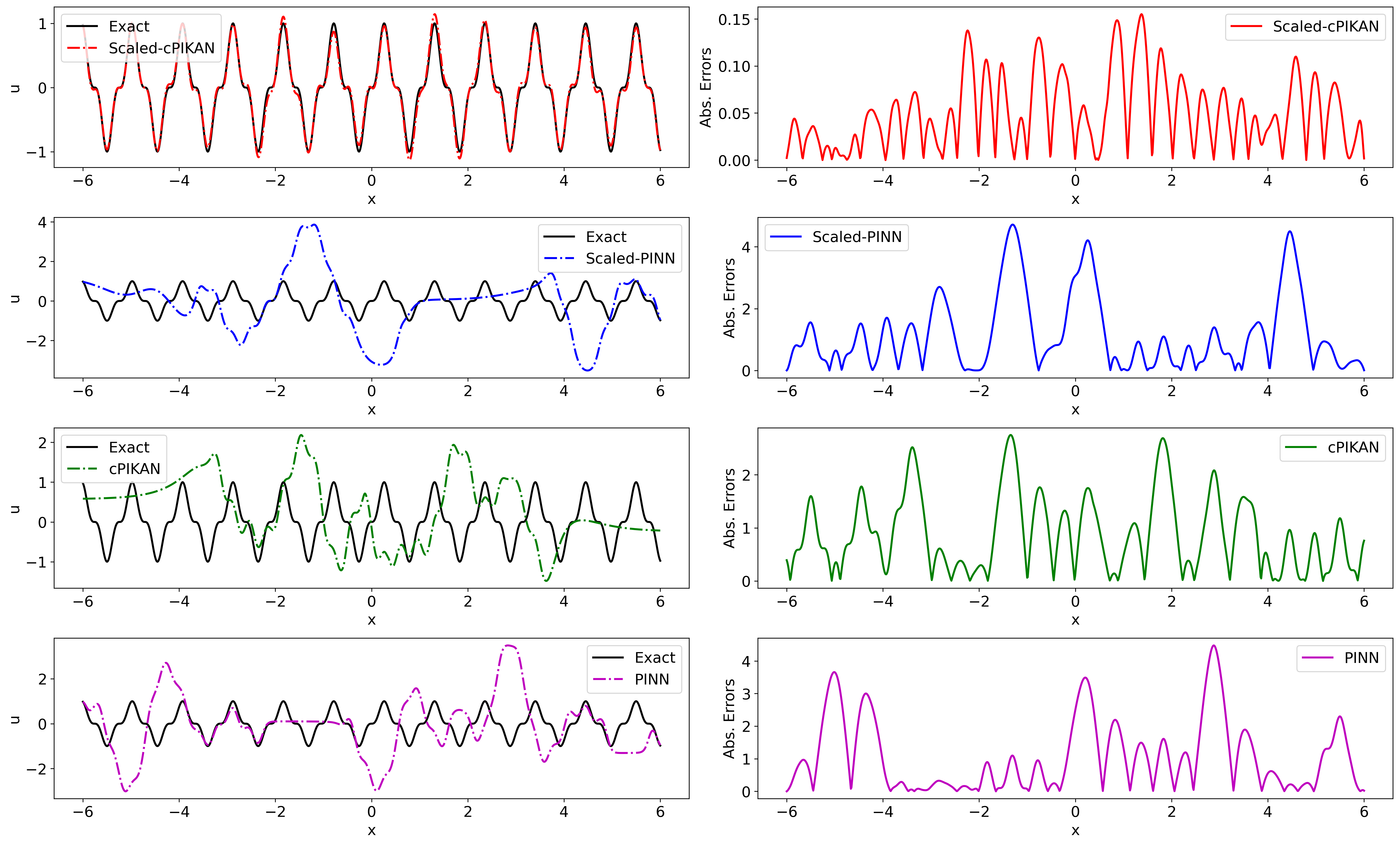}
    \caption{Inverse problem of the reaction-diffusion equation (Example \ref{Exam.RD}) results for \(u\) when \(\delta_u = 0.05\) and \(\delta_f = 0.05\) using Scaled-cPIKAN, Scaled-PINN, cPIKAN, and PINN methods (from top to bottom). Left Panels: comparison between the ground truth solution and predictions of the methods. Right Panels: absolute error distributions.}
    \label{fig:INV_RD_N2}
\end{figure}
%%%=============================================================
% % % ===========================================
\begin{figure}[!h]
    \centering
    \includegraphics[width=0.95\linewidth]{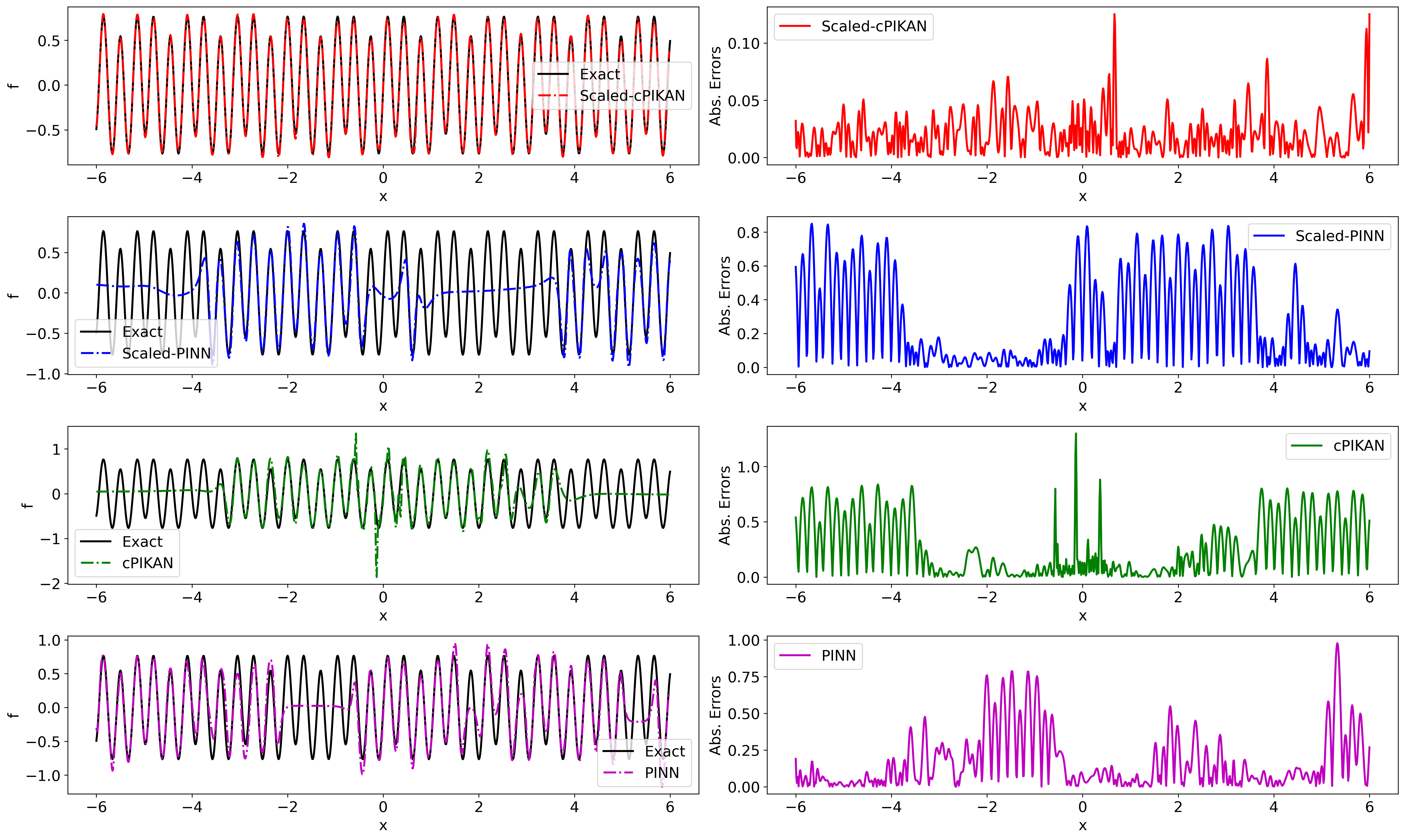}
    \caption{Inverse problem of the reaction-diffusion equation (Example \ref{Exam.RD}) results for \(f\) when \(\delta_u = 0.05\) and \(\delta_f = 0.05\) using Scaled-cPIKAN, Scaled-PINN, cPIKAN, and PINN methods (from top to bottom). Left Panels: comparison between the ground truth solution and predictions of the methods. Right Panels: absolute error distributions.}
    \label{fig:INVf_RD_N2}
\end{figure}
% % % ===========================================

Figure \ref{fig:INV_RD_Loss} provides the training loss history of the inverse problem with different noise levels in the inverse version of the reaction-diffusion equation.
The left panel corresponds to the noise-free case (\(\delta_u = 0.00\), \(\delta_f = 0.00\)), where Scaled-cPIKAN has the fastest convergence and reaches the lowest loss among all methods. More precisely, the loss value of Scaled-cPIKAN goes below \(10^{-4}\), while the other three methods stabilize around \(10^{-1}\), highlighting the superior efficiency and precision of Scaled-cPIKAN in this setting.
The middle panel depicts the setting of noisy labeled data (\(\delta_u = 0.05\), \(\delta_f = 0.00\)), where Scaled-cPIKAN is still the most accurate method. In this case, Scaled-cPIKAN reduces the loss to \(10^{-4}\), while the vanilla cPIKAN struggles to drop below \(10^{-1}\), and the other PINN-based methods achieve losses only as low as \(10^{-2}\).
The right panel depicts the most challenging case, with noise present in both the supervised data and the source term (\(\delta_u = 0.05\), \(\delta_f = 0.05\)). Yet, Scaled-cPIKAN demonstrates remarkable robustness and maintains its superior convergence, reducing the loss below \(10^{-3}\), while Scaled-PINN and cPIKAN stabilize around \(10^{-1}\), and PINN stabilizes around \(10^{-2}\). This underlines the robustness and the accuracy of Scaled-cPIKAN for solving the inverse problem of the reaction-diffusion equation even with noisy conditions.
% % % 

Figures \ref{fig:INV_RD_N2} and \ref{fig:INVf_RD_N2} present the inverse problem results for the reaction-diffusion equation using the Scaled-cPIKAN, Scaled-PINN, cPIKAN, and PINN methods. The plots represent the most challenging case with noise in both the solution and the source term, that is, \((\delta_u, \delta_f) = (0.05, 0.05)\).
In Fig.~\ref{fig:INV_RD_N2}, among all methods, the results of Scaled-cPIKAN are the closest to the ground truth solution, taking into account the noisy input provided. The robustness of Scaled-cPIKAN shows clearly, being able to stay well below $0.15$ for the maximum absolute error of approximation of \(u\) - significantly outperforming methods that all have maximum errors above $4.0$.
Similarly, in Fig.~\ref{fig:INVf_RD_N2}, Scaled-cPIKAN obtains a maximum absolute error in approximating \(f\) of less than $0.12$ while the corresponding errors for the other methods exceed $1.0$. 
These results highlight the great stability and accuracy of the Scaled-cPIKAN framework when solving the noisy inverse reaction-diffusion problem. Meanwhile, Scaled-PINN and PINN are unstable, and they cannot provide good approximations for higher noise. Besides, vanilla cPIKAN's performance degrades close to the boundaries, which further indicates the advantages of the scaled version.
For the approximation of \( u \) when \((\delta_u, \delta_f) = (0.00, 0.00)\) and \((\delta_u, \delta_f) = (0.05, 0.00)\), refer to Figs.~\ref{fig:INV_RD_N0} and \ref{fig:INV_RD_N1} in \ref{App.Fig}. Similarly, for the approximation of \( f \) in the same cases, see Figs.~\ref{fig:INVf_RD_N0} and \ref{fig:INVf_RD_N1} in  \ref{App.Fig}.
\begin{figure}[h!]
    \centering
    \includegraphics[width=0.8\linewidth]{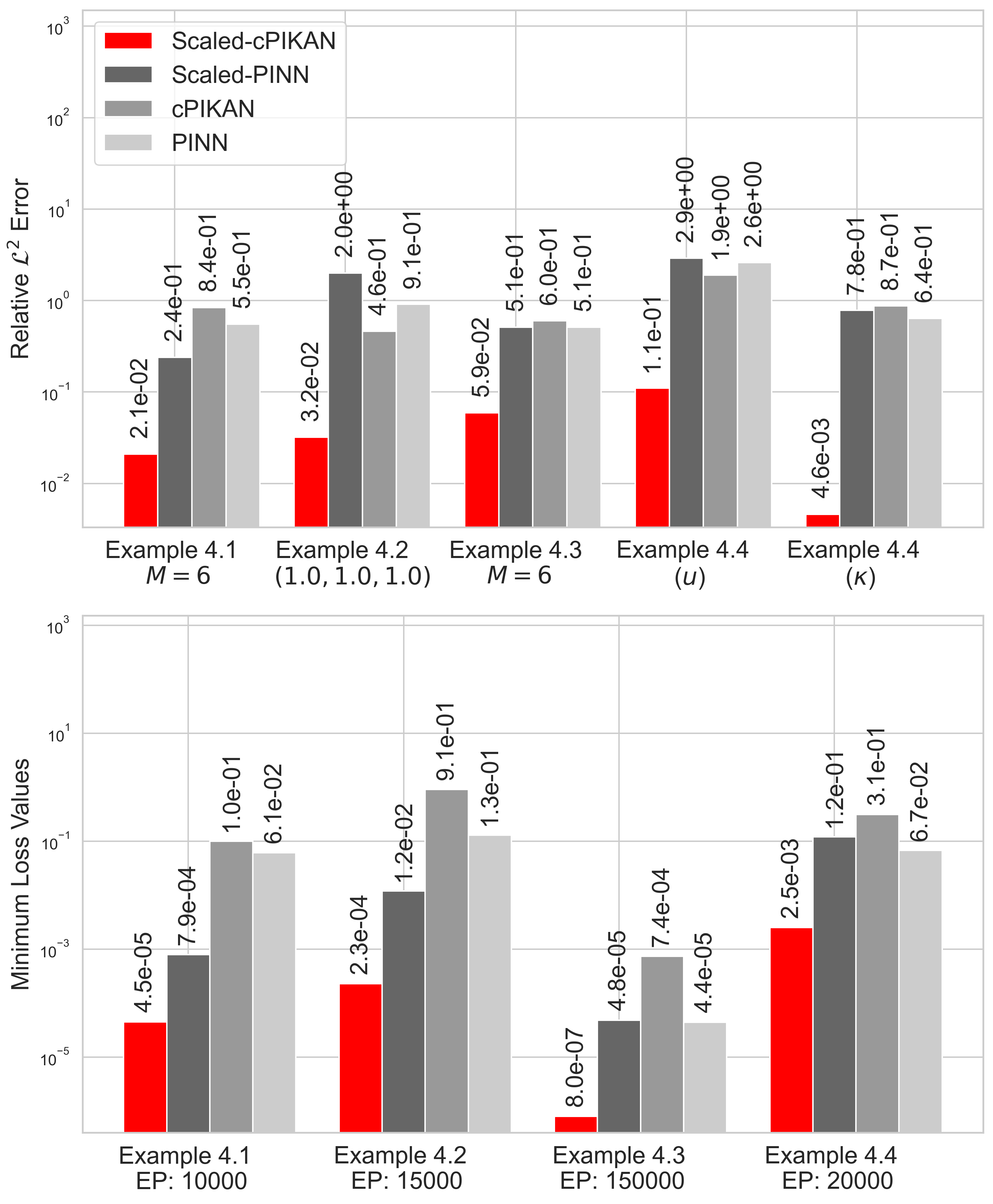}
    \caption{
    Overall comparison of the performance of Scaled-cPIKAN, Scaled-PINN, cPIKAN, and PINN on challenging examples reported in Sec.~\ref{Sec.res}. The top panel compares the relative \(\mathcal{L}^2\) error, while the bottom panel highlights the minimum loss values at the specified number of epochs (EP), showcasing the superior accuracy and lower loss values of Scaled-cPIKAN in all cases.
    % Example 4.1 considers an extended spatial domain with \(M = 6\). Example 4.2 examines the settings \((a_1, a_2, \kappa) = (1.0, 1.0, 1.0)\). Example 4.3 represents the unsupervised version on an extended domain with \(M = 6\). Example 4.4 focuses on the inverse problem with noisy data \((\delta_u, \delta_f) = (0.05, 0.05)\), addressing the approximation of the solution \(u\) and PDE parameter \(\kappa\). 
    }
    \label{fig:Summary}
\end{figure}

In this work, the performance of Scaled-cPIKAN, Scaled-PINN, cPIKAN, and PINN is compared for solving PDEs such as the diffusion equation, Helmholtz equation, Allen-Cahn equation, and reaction-diffusion equation in both forward and inverse settings. For each problem, the most complicated and challenging configurations have been selected to test the accuracy and robustness of the proposed methods. Figure~\ref{fig:Summary} summarizes the results, showing that Scaled-cPIKAN outperforms all other methods by a wide margin in terms of lower relative $\mathcal{L}^2$ errors and minimum loss values obtained. These results confirm that Scaled-cPIKAN is an efficient and robust solver for partial differential equations on extended spatial domains and, in general, when oscillatory behaviors and other challenging dynamics are present, which are difficult to resolve with conventional methods. Note that we mainly considered the performance of Scaled-cPIKAN within spatial domains up to \([-6, 6]\), where the method demonstrated excellent accuracy. We observed changes in the solution dynamics for domains considerably larger (for example, extending beyond 100) that provide many opportunities for further refinement and optimization of the approach. These observations open up investigations for strategies such as domain decomposition, whereby larger domains could be decomposed into smaller sub-domains mapped to \([-1, 1]\), such that further decomposition optimally balances computational efficiency with accuracy. Moreover, the focus of this study is on extended spatial domains, the temporal domain is kept within \([0, 1]\). The extension of our scaling method for larger temporal intervals or to higher-dimensional PDEs is an exciting avenue for future research. One very promising direction is the incorporation of adaptive scaling mechanisms that can dynamically adjust to different solution complexities across both spatial and temporal domains, thus ensuring robust and efficient performance under a wide range of conditions.

%%%=============================================================
\section{Conclusion}\label{Sec.Conclusion}

In this paper, we introduced Scaled-cPIKAN, a novel neural network framework that integrates domain scaling into the cPIKAN architecture.
Scaled-cPIKAN transforms spatial variables in the governing equations to a standardized domain \([-1,1]^d\) as required by Chebyshev  polynomials, allowing the framework to efficiently represent oscillatory dynamics over extended spatial domains while improving computational performance.
Our motivation for proposing Scaled-cPIKAN stemmed from the limitations of existing approaches, such as Scaled-PINN, cPIKAN, and PINN, which often struggle with accuracy and convergence when applied to large domains, especially in inverse problems or scenarios involving noisy data.
For validation purposes, we evaluated Scaled-cPIKAN across four standard benchmark problems: the diffusion equation, the Helmholtz equation, the Allen-Cahn equation, as well as both forward and inverse configurations of the reaction-diffusion equation, under conditions both with and without noise. These tests focused on assessing the accuracy of the solution and the convergence behavior of the loss function in a range of diverse scenarios.
Scaled-cPIKAN consistently outperformed the other methods, delivering several orders of magnitude higher accuracy and faster convergence. For example, in a challenging inverse reaction-diffusion problem with noisy data, Scaled-cPIKAN reduced errors in parameter estimation and solution approximation by an average of 97\% compared to Scaled-PINN, cPIKAN, and PINN. It also achieved loss values below \(10^{-4}\), while other methods plateaued around
\(10^{-1}\). 
In the future, we hope to improve the scalability of Scaled-cPIKAN so that it can tackle PDEs on much larger
spatial and temporal domains, broadening its potential applications. 

\section{Acknowledgements}
S.A.F. acknowledges the support by the U.S. Department of Energy’s Office of Environmental Management  (award no.: DE-EM0005314).

\section{Conflict of Interest}
The authors declare no conflict of interests.

%PPPPPPPPPPPPPPPPPPPPPPPPPPPPPPPPPPPPPPPPPPPPPPPPPPPPPPPPPPPPPPPPPPPPPPPPPPPPPPPPPPPPPPPPPPPPPPPPPPPPPPPPPPPPPPPP
\def\mybibdoicolor{\color{black}}
\newcommand*{\doi}[1]{\href{\detokenize{#1}} {\raggedright\mybibdoicolor{DOI: \detokenize{#1}}}}

\bibliographystyle{unsrtnat}
\bibliography{references.bib}

%%%%%%========================================================
\newpage
\appendix
\renewcommand{\thefigure}{A.\arabic{figure}}
\setcounter{figure}{0}

\section{Additional Supporting Figures}\label{App.Fig}
%%%%%%========================================================
\begin{figure}[!h]
    \centering
    \includegraphics[width=1.0\linewidth]{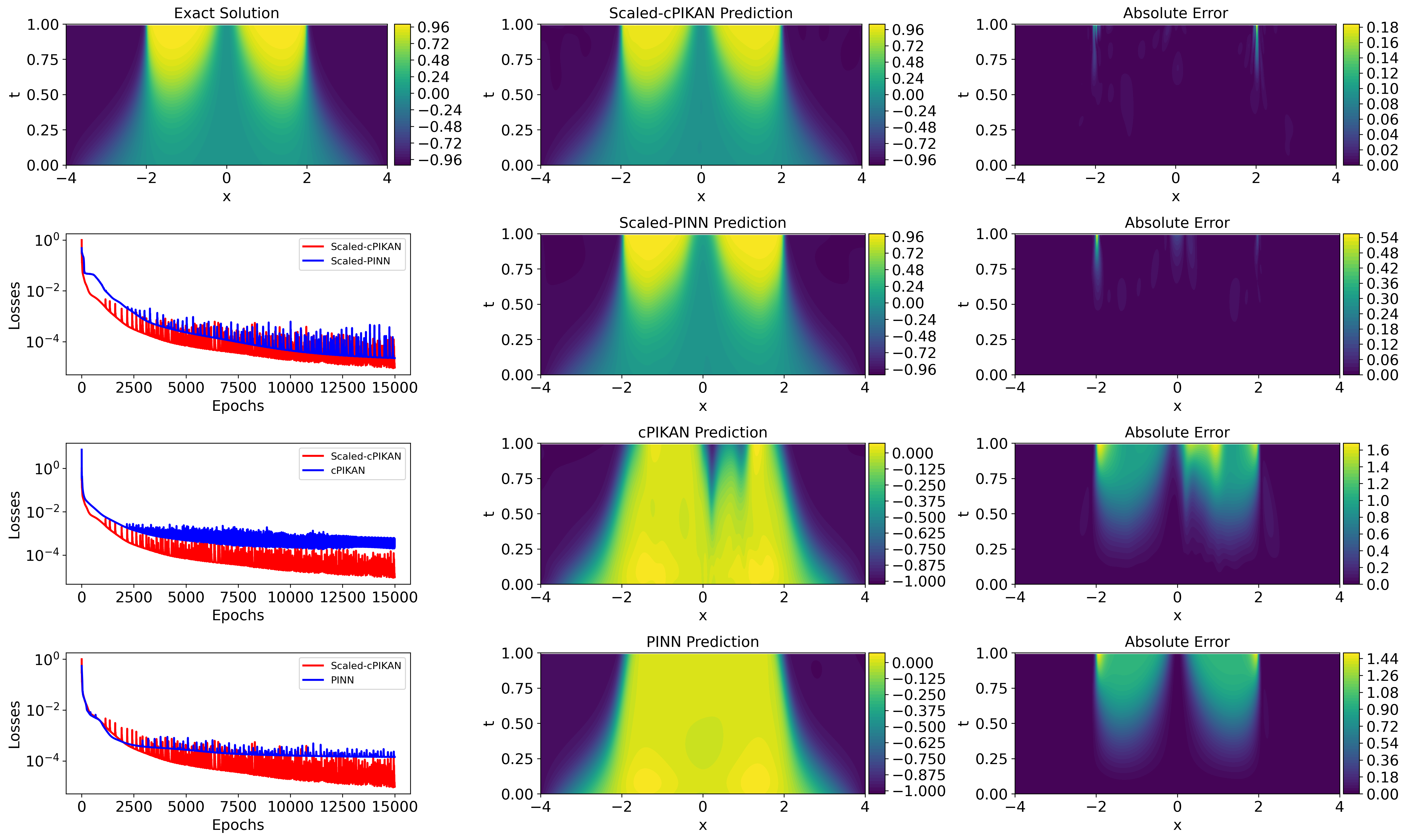}
    \caption{Comparison of the solutions predicted for the Allen-Cahn equation (Example \ref{Exam.AC}) in a supervised version on \([-4, 4]\times [0,1]\). Each figure displays the outcomes of Scaled-cPIKAN, Scaled-PINN, cPIKAN, and PINN (from top to bottom). From left to right in each row: the loss function, the predicted solution \(u\), and the absolute error of the prediction are shown. In the first row, the first plot shows the ground truth solution instead of the loss function for comparison.}
    \label{fig:AC4}
\end{figure}
%%%%%%========================================================
%%%%%%========================================================
\begin{figure}
    \centering
    \includegraphics[width=1.0\linewidth]{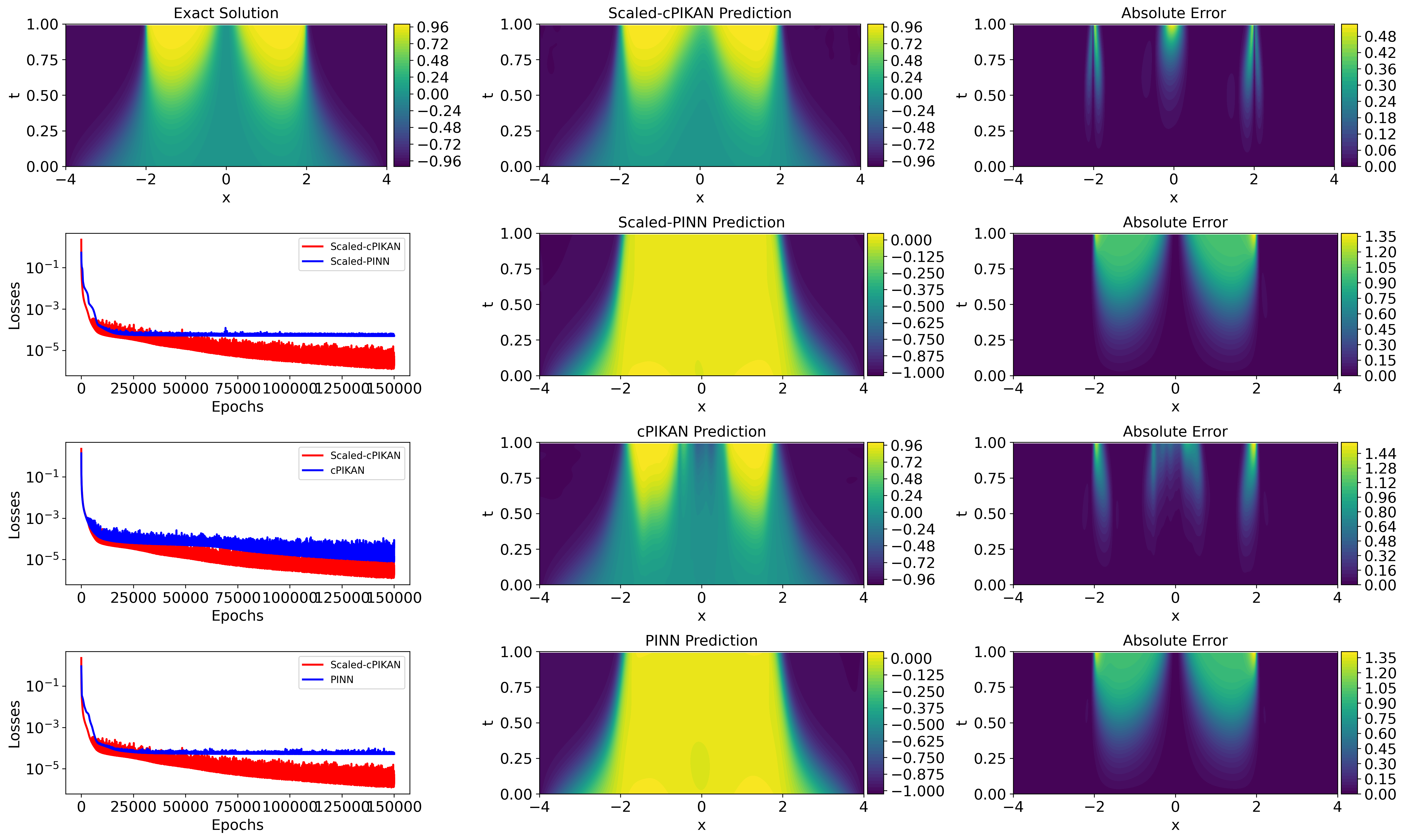}
    \caption{Comparison of the solutions predicted for the Allen-Cahn equation (Example \ref{Exam.AC}) in a unsupervised version on \([-4, 4]\times [0,1]\). Each figure displays the outcomes of Scaled-cPIKAN, Scaled-PINN, cPIKAN, and PINN (from top to bottom). From left to right in each row: the loss function, the predicted solution \(u\), and the absolute error of the prediction are shown. In the first row, the first plot shows the ground truth solution instead of the loss function for comparison.}
    \label{fig:unsAC4}
\end{figure}
%%%%%%========================================================
% % % ===========================================
\begin{figure}
    \centering
    \includegraphics[width=1.0\linewidth]{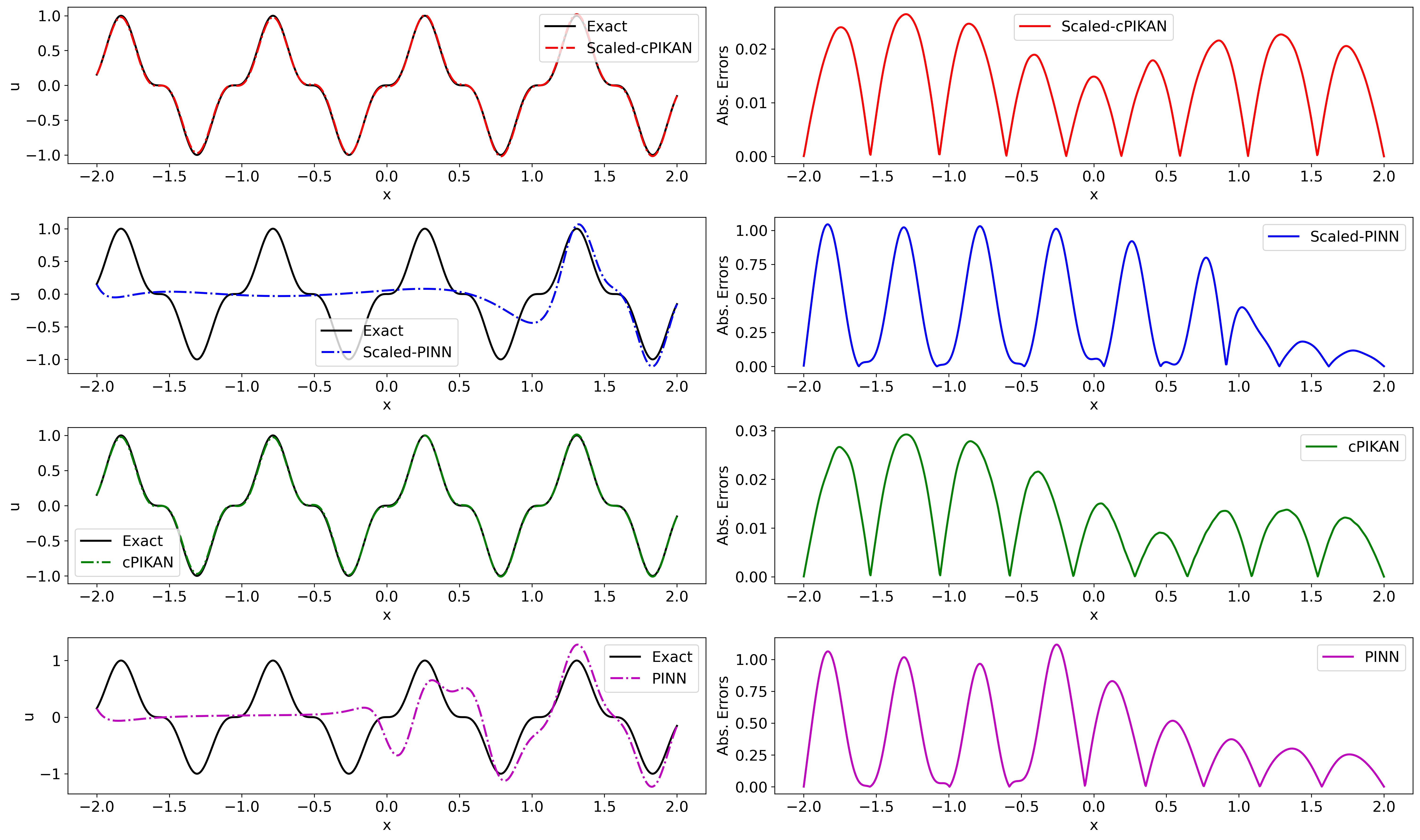}
    \caption{Forward solutions of the reaction-diffusion equation (Example \ref{Exam.RD}) for $M=2$ using Scaled-cPIKAN, Scaled-PINN, cPIKAN, and PINN methods (from top to bottom). Left Panels: comparison between the ground truth solution and predictions of the methods. Right Panels: absolute error distributions.}
    \label{fig:FRD2}
\end{figure}
\begin{figure}
    \centering
    \includegraphics[width=1.0\linewidth]{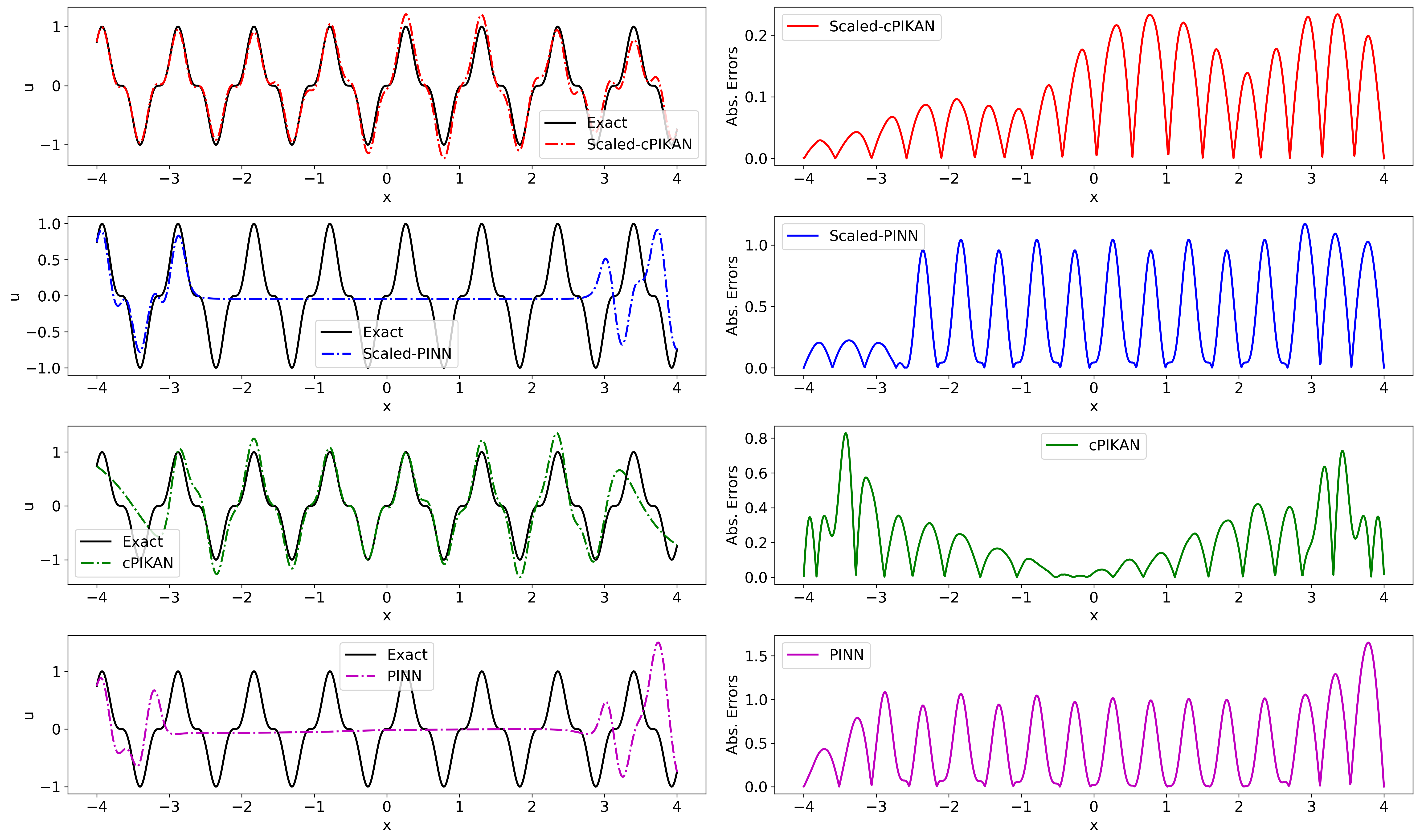}
    \caption{Forward solutions of the reaction-diffusion equation (Example \ref{Exam.RD}) for $M=4$ using Scaled-cPIKAN, Scaled-PINN, cPIKAN, and PINN methods (from top to bottom). Left Panels: comparison between the ground truth solution and predictions of the methods. Right Panels: absolute error distributions.}
    \label{fig:FRD4}
\end{figure}
%%%%%%========================================================
% % % ===========================================
\begin{figure}
    \centering
    \includegraphics[width=1.0\linewidth]{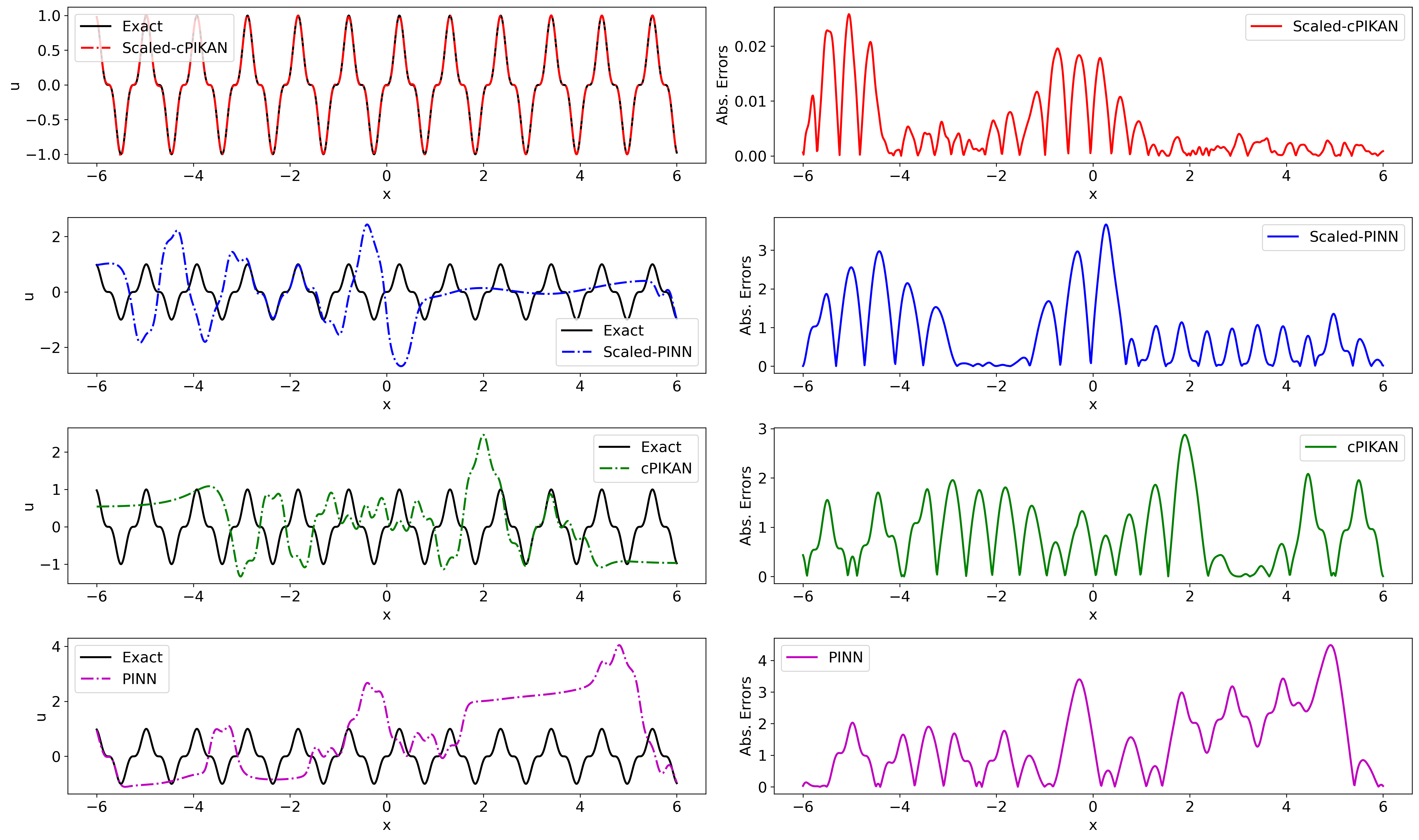}
    \caption{(Inverse problem of the reaction-diffusion equation (Example \ref{Exam.RD}) results for \(u\) when \(\delta_u = 0.0\) and \(\delta_f = 0.0\) using Scaled-cPIKAN, Scaled-PINN, cPIKAN, and PINN methods (from top to bottom). Left Panels: comparison between the ground truth solution and predictions of the methods. Right Panels: absolute error distributions.}
    \label{fig:INV_RD_N0}
\end{figure}
%%%=============================================================
% % % ===========================================
\begin{figure}
    \centering
    \includegraphics[width=1.0\linewidth]{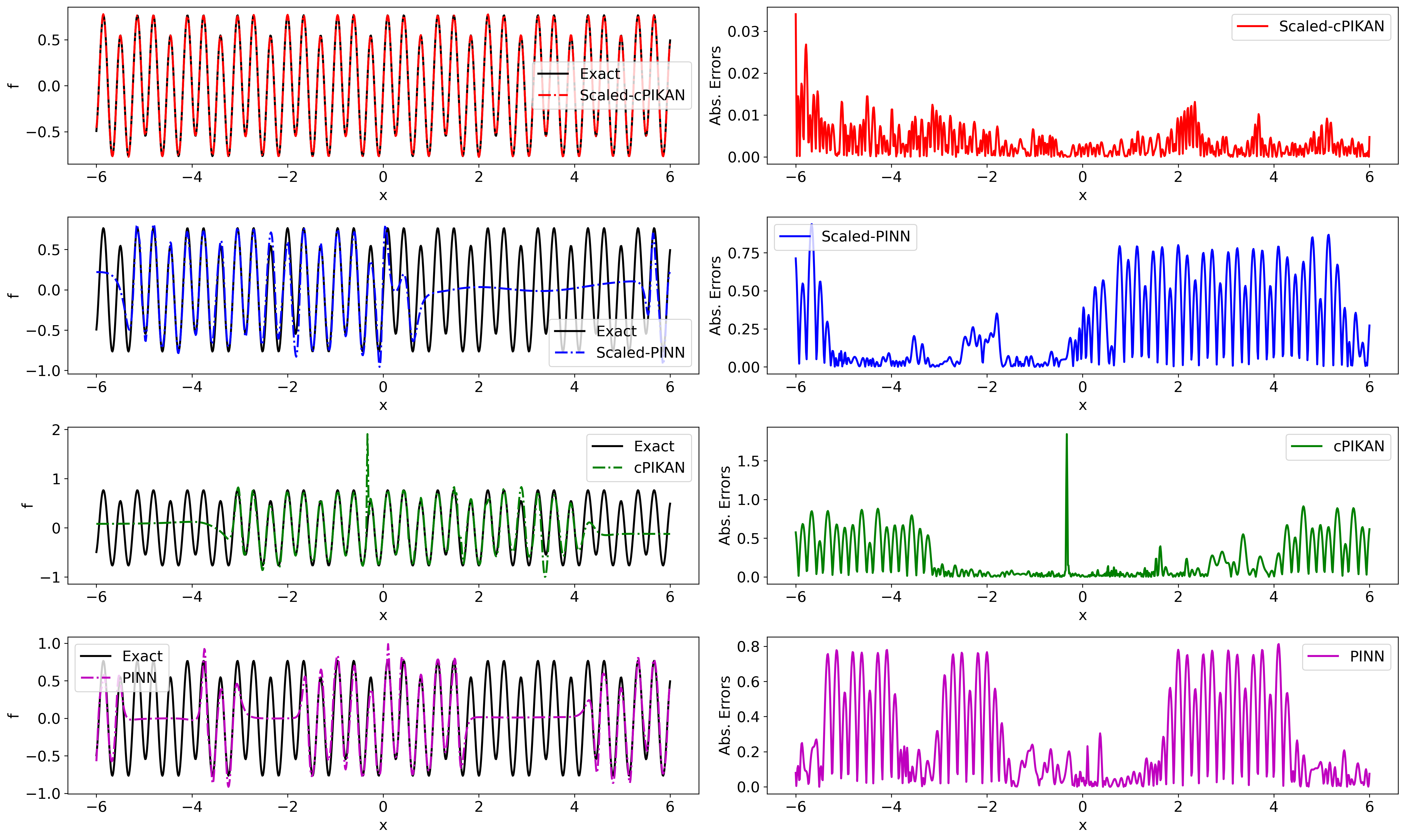}
    \caption{Inverse problem of the reaction-diffusion equation (Example \ref{Exam.RD}) results for \(f\) when \(\delta_u = 0.00\) and \(\delta_f = 0.00\) using Scaled-cPIKAN, Scaled-PINN, cPIKAN, and PINN methods (from top to bottom). Left Panels: comparison between the ground truth solution and predictions of the methods. Right Panels: absolute error distributions.}
    \label{fig:INVf_RD_N0}
\end{figure}
%%%=============================================================

% % % ===========================================
\begin{figure}
    \centering
    \includegraphics[width=1.0\linewidth]{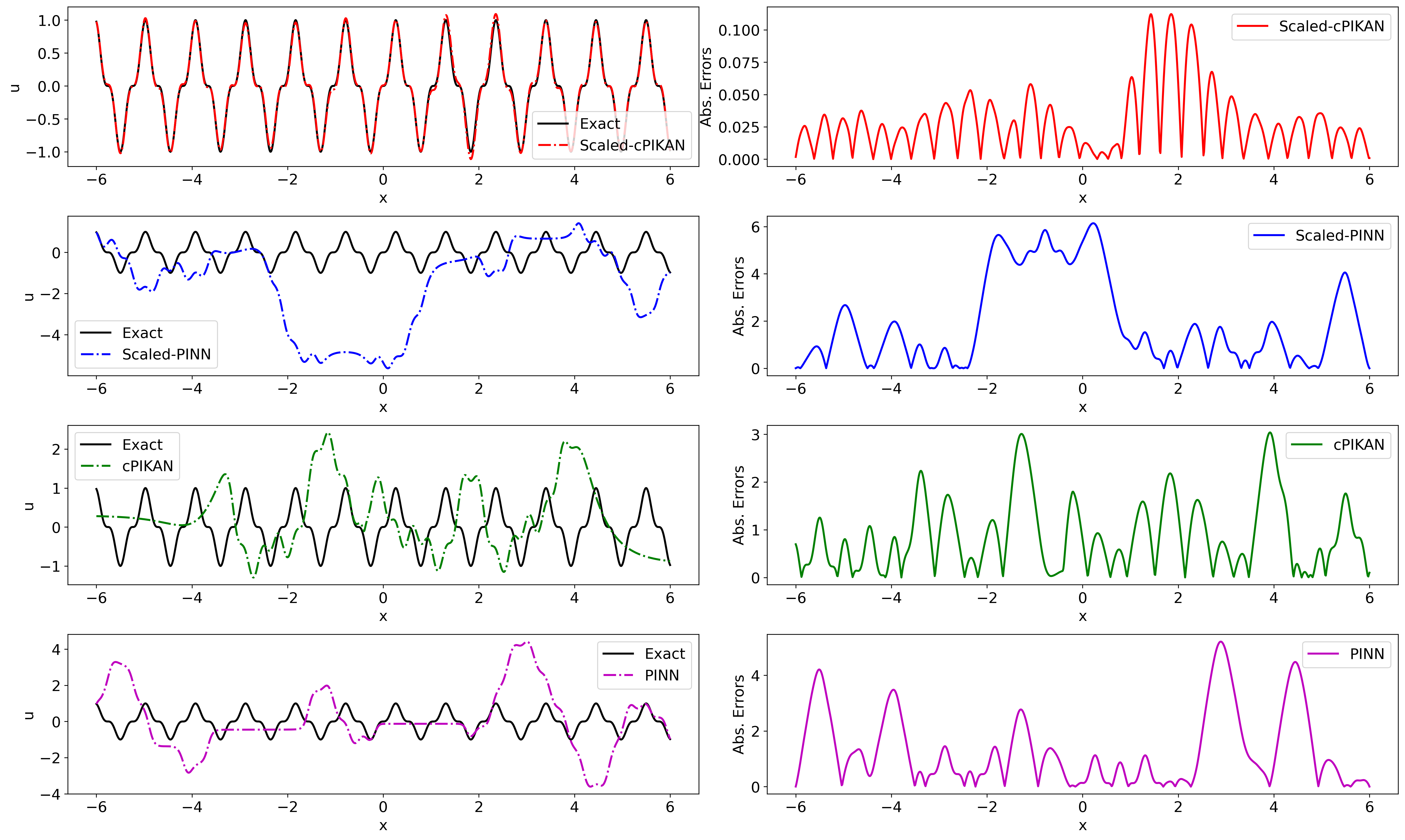}
    \caption{Inverse problem of the reaction-diffusion equation (Example \ref{Exam.RD}) results for \(u\) when \(\delta_u = 0.05\) and \(\delta_f = 0.00\) using Scaled-cPIKAN, Scaled-PINN, cPIKAN, and PINN methods (from top to bottom). Left Panels: comparison between the ground truth solution and predictions of the methods. Right Panels: absolute error distributions.}
    \label{fig:INV_RD_N1}
\end{figure}
%%%=============================================================
% % % ===========================================
\begin{figure}
    \centering
    \includegraphics[width=1.0\linewidth]{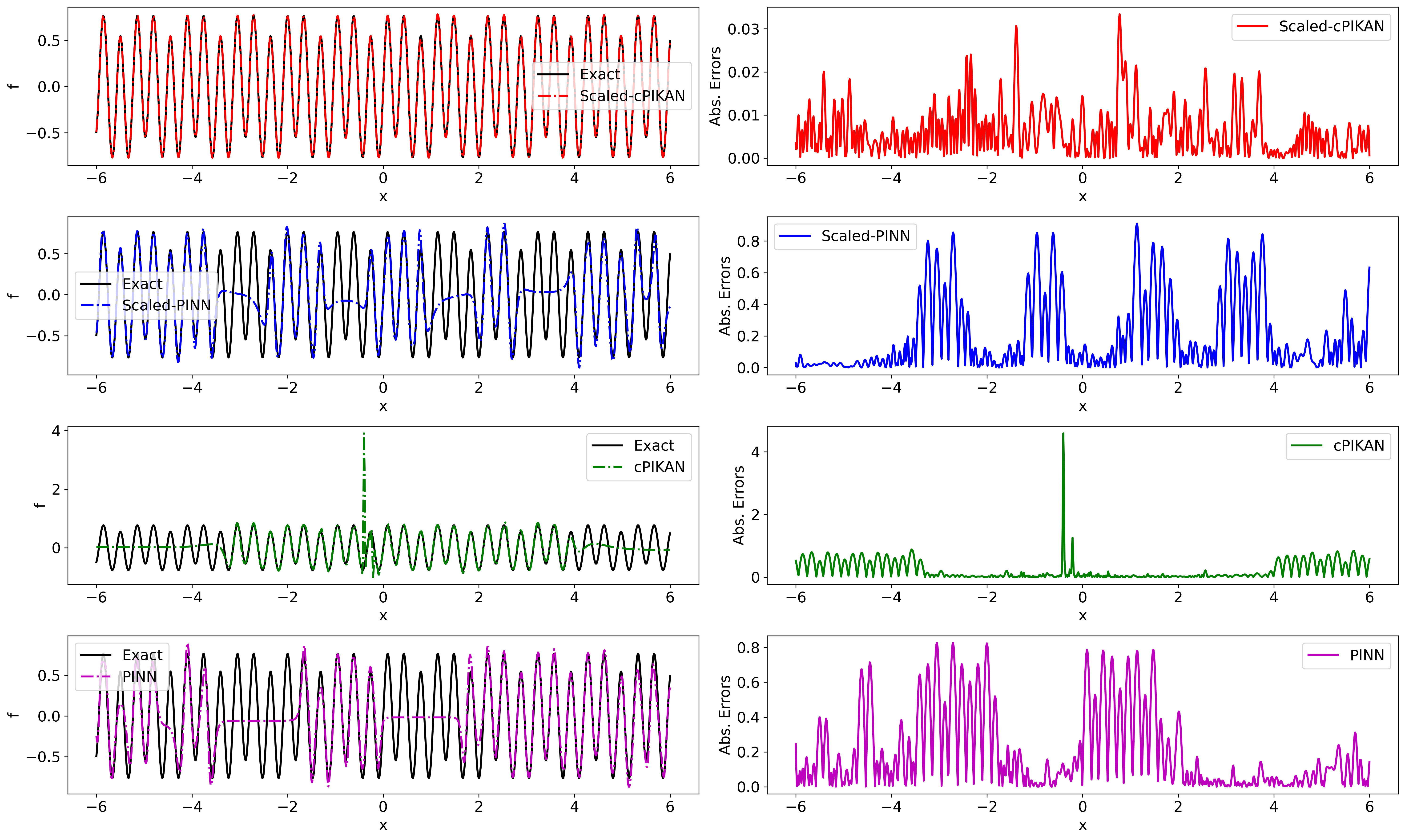}
    \caption{Inverse problem of the reaction-diffusion equation (Example \ref{Exam.RD}) results for \(f\) when \(\delta_u = 0.05\) and \(\delta_f = 0.00\) using Scaled-cPIKAN, Scaled-PINN, cPIKAN, and PINN methods (from top to bottom). Left Panels: comparison between the ground truth solution and predictions of the methods. Right Panels: absolute error distributions.}
    \label{fig:INVf_RD_N1}
\end{figure}
%%%=============================================================

\end{document}